\documentclass{amsart}

\usepackage{amsfonts}
\usepackage{amsmath}
\usepackage{amsthm}
\usepackage{amssymb}
\usepackage[english]{babel}
\usepackage{graphics}
\usepackage{latexsym}
\usepackage{longtable} \usepackage{mathrsfs}
\usepackage{graphicx}
\usepackage{wrapfig} 
\usepackage[pdftex,colorlinks=true]{hyperref}

\usepackage{fullpage}
\usepackage{tristan}
\usepackage{subfig}

\newcommand{\mD}{\mathcal{D}}

\renewcommand{\S}{\mathcal{S}}

\renewcommand{\R}{\mathbb{R}}
\renewcommand{\N}{\mathbb{N}}

\renewcommand{\X}{\textbf{X}}

\newcommand{\noi}{\noindent}
\renewcommand{\ms}{\medskip}

\newcommand{\al}{\alpha}
\newcommand{\be}{\beta}

\renewcommand{\de}{\delta}
\renewcommand{\De}{\Delta}

\newcommand{\Om}{\Omega}

\newcommand{\weak }{\, -\!\!\!\!-\!\!\!\!\rightharpoonup}
\newcommand{\weakstar }{ \overset{\, *_{\phantom{|}}}{{\smash{\weak }}\, } }

\newcommand{\larrow }{\longrightarrow}
\newcommand{\ot}{\otimes}

\newcommand{\lmapsto}{\longmapsto}
\newcommand{\ri}{\rightarrow}
\renewcommand{\p}{\partial}
\newcommand{\sub}{\subseteq}
\newcommand{\set}{\setminus}
\newcommand{\by}{\times}
\newcommand{\rk}{\textrm{rk}}

\newcommand{\inter}{\textrm{int}}

\newcommand{\bt}{\begin{theorem}}\newcommand{\et}{\end{theorem}}
\newcommand{\bd}{\begin{definition}}\newcommand{\ed}{\end{definition}}
\newcommand{\bl}{\begin{lemma}}\newcommand{\el}{\end{lemma}}
\newcommand{\beq}{\begin{equation}}\newcommand{\eeq}{\end{equation}}
\newcommand{\bc}{\begin{claim}}\newcommand{\ec}{\end{claim}}
\newcommand{\bex}{\begin{example}}\newcommand{\eex}{\end{example}}
\newcommand{\bcor}{\begin{corollary}}\newcommand{\ecor}{\end{corollary}}
\newcommand{\bp}{\begin{proof}}\renewcommand{\ep}{\end{proof}}

\numberwithin{equation}{section}

\setboolean{shownotes}{true}

\begin{document}

\title{On the Numerical Approximation of $\infty$-Harmonic Mappings}

\author{Nikos Katzourakis and Tristan Pryer}
\address{Department of Mathematics and Statistics, University of Reading, Whiteknights, PO Box 220, Reading RG6 6AX, United Kingdom.}
\email{n.katzourakis@reading.ac.uk, t.pryer@reading.ac.uk}


\date{}


\begin{abstract} Given a map $u : \Om \sub \R^n \larrow  \R^N$, the $\infty$-Laplacian is the system
\[  \label{1}
\De_\infty u  \, :=\, \Big(\text{D}u \ot \text{D}u  + |\text{D}u|^2 [\text{D}u]^\bot \! \ot I \Big) : \text{D}^2 u\, = \, 0.  \tag{1}
\]
\eqref{1} is the model system of vectorial Calculus of Variations in $L^\infty$ and arises as the ``Euler-Lagrange'' equation in relation to the supremal functional 
\[  \label{2}
E_\infty(u,\Om)\, :=\,  \| \text{D}u \|_{L^\infty(\Om)}.    \tag{2}
\]
The scalar case of \eqref{1} has been introduced by Aronsson in the 1960s and by now is relatively classical and well understood. The general system \eqref{1} has been discovered and studied by the first author in a series of recent papers. Supremal functionals are fundamental for applications because they provide more realistic models as opposed to conventional integral models. Herein we provide numerical approximations of solutions to the Dirichlet problem when $n=2$ and $N=2,3$ for certain carefully selected boundary data on the unit square. Our experiments demonstrate interesting and unexpected phenomena occurring in $L^\infty$ and provide insights on the structure of general solutions and the natural separation to phases they present.
\end{abstract}

\maketitle

\section{Introduction} \label{section1}

Let $n,N\in \N$ and $\Om$ an open set in $\R^n$. Given a smooth map $u:\Om \sub \R^n \larrow  \R^N$ with components $(u_1,...,u_N)^\top$, the gradient matrix map $\text{D}u :\Om \sub \R^n \larrow  \R^{Nn}$ is denoted by $\text{D}u=(\text{D}_iu_\al)^{\al=1...N}_{i=1...n}$ and the hessian tensor $\text{D}^2u : \Om \sub \R^n \larrow  \R^{Nn^2}_s$ is denoted by $\text{D}^2u=(\text{D}^2_{i j}u_\al)^{\al=1...N}_{i,j=1...n}$. We use $\R^{Nn}$ and $\R^{Nn^2}_s$ to denote respectively the matrix space and the symmetric tensor space in which $\text{D}u$ and $\text{D}^2u$ are valued. In this paper we are interested in the numerical approximation of solutions to the $\infty$-Laplacian which is defined on smooth maps as the following PDE system
\beq \label{1.1}
\De_\infty u  \, :=\, \Big(\text{D}u \ot \text{D}u  + |\text{D}u|^2 [\text{D}u]^\bot \! \ot I \Big) : \text{D}^2 u\, = \, 0.
\eeq
In the above, ``$\ot$" denotes the usual tensor product, ``$|\text{D}u|$" denotes the Euclidean norm of the gradient in the matrix space $\R^{Nn}$, that is 
\[
|\text{D}u|^2\, =\, \sum_{\al=1}^N\sum_{i=1}^n \text{D}_iu_\al \, \text{D}_iu_\al\,  \equiv \, \text{D}u:\text{D}u, 
\]
``$I$" is the identity matrix and ``$[\text{D}u]^\bot$" denotes the orthogonal projection on the subspace of $\R^N$ which is the orthogonal complement of the range of the linear mapping $\text{D}u(x) :\R^n \larrow  \R^N$ for fixed $x\in \Om$:
\beq \label{1.2}
[\text{D}u]^\bot\, :=\, \text{Proj}_{(R(\text{D}u))^\bot}.
\eeq
The notation ``$:$" symbolises a contraction with respect to 3 indices as above which extends the usual Euclidean inner product in $\R^{Nn}$.  In index form with respect to the standard bases, \eqref{1.1} reads
\[
\sum_{\be=1}^N\sum_{i,j=1}^n\Big( \text{D}_i u_\al \, \text{D}_ju_\be \,+ \, |\text{D}u|^2[\text{D}u]^\bot_{\al \be} \,\de_{ij}\Big) \text{D}^2_{ij}u_\be\, =\, 0, \ \ \ \ \al=1,...,N.
\]
The system \eqref{1.1} (whose solution we call $\infty$-Harmonic mappings) arises as the analogue of the Euler-Lagrange equations associated to the functional
\beq \label{1.3}
E_\infty(u,\Om)\, :=\,  \big\| |\text{D}u| \big\|_{L^\infty(\Om)}.   
\eeq
The so-called supremal functional \eqref{1.3} and the PDE system \eqref{1.1} are the archetypal model objects of vectorial Calculus of Variations in the space $L^\infty$. In the scalar case $N=1$, the system \eqref{1.1} simplifies to the following single equation
\beq \label{1.4}
\text{D}u\ot \text{D}u :\text{D}^2u\, = \sum_{i,j=1}^n\text{D}_i u\, \text{D}_j u \,\text{D}^2_{ij}u\, =\,0
\eeq
since the second term involving the orthogonal projection \eqref{1.2} vanishes identically.  The field of Calculus of Variations in the space $L^\infty$ has been pioneered by Aronsson in the 1960s \cite{A1}-\cite{A6} who was the first to derive and study \eqref{1.4} is relation to variational problems arising from \eqref{1.3}. Since then, the field has undergone an explosion of interest within the PDE community, we refer e.g.\ to \cite{B,C,K} and references therein. Until the early 2010s, all considerations where restricted to the scalar case of $N=1$. The general vectorial case of \eqref{1.1} as well as the study of the associated PDE systems arising from more general first order functionals
\beq 
 \label{1.5}
 E_\infty(u,\Om)\, :=\,  \big\| H(\cdot,u,\text{D}u) \big\|_{L^\infty(\Om)} 
\eeq
has been pioneered by the first author in a series of recent papers \cite{K1}-\cite{K9}. Except for the intrinsic mathematical interest of the field, $L^\infty$ functionals are very important for applications as well since we obtain more realistic variational models as opposed to their integral counterparts: for example in mechanics it is preferable to know to maximum tolerance of a specific material to distractive forces rather than the average tolerance.

A basic difficulty arising already in the scalar case is that \eqref{1.4} is degenerate elliptic and in non-divergence form and classical approaches to define and study weak solutions via integration-by-parts fail. In general the Dirichlet problem can not be solved in the class of smooth functions since Aronsson himself demonstrated the existence of singular solutions in \cite{A6,A7} which however are minimising for the functional. The theory of viscosity solutions of Crandall-Ishii-Lions (see \cite{CIL, C} and for a pedagogical introduction we refer to \cite{K}) proved to be the appropriate framework to study scalar $L^\infty$ variational problems. 

In the vectorial case, $N\geq 2$, singular solutions of \eqref{1.1} still appear (see \cite{K1,K3}). A further difficulty associated to \eqref{1.1} which is not present when $N=1$ is that the projection $[\text{D}u]^\bot$ may be discontinuous even for $C^\infty$ maps $u$, whence the nonlinear operator $\De_\infty$ of \eqref{1.2} may have discontinuous coefficients even when applied to $C^\infty$ maps. This happens because the range of the gradient matrix $\text{D}u(x)\in \R^{Nn}$ may not be constant throughout the domain $\Om$. A simple example of $C^\infty$ smooth solution to \eqref{1.1} from the unit rectangle of $\R^2$ into $\R^2$ is $u(x,y) = e^{ix}-e^{iy}$ (\cite{K1}). The rank of $\text{D}u$ is $2$ away from the diagonal $\{x=y\}$ because the partial derivative $\text{D}_x u$, $\text{D}_y u$ are linearly independent but on the diagonal the rank jumps to $1$ and \eqref{1.2} becomes discontinuous thereon. This is a general phenomenon and the solutions of \eqref{1.1} always come in phases of different rank values of the gradient with interfaces separating them (see \cite{K2,K3}). The appropriate duality-free notion of generalised solution for \eqref{1.1} has very recently been proposed by the first author in \cite{K7} and is based on the probabilistic representation of the hessian which does not exist classically via Young measures valued the compactification of the space $\R^{Nn^2}_s$. In this setting, among other far-reaching results, we have been able to prove existence to the Dirichlet problem for \eqref{1.1} and for the equations arising from \eqref{1.5} for $n=1$. Subsequent developments in the context of this theory appear in \cite{K8}-\cite{K11}.

In the scalar case some numerical schemes have been proposed for the direct approximation of viscosity solutions of \eqref{1.4}. In \cite{Oberman:2005,EsedogluOberman:2011,Oberman:2013} Oberman uses techniques from Barles and Souganidis \cite{BarlesSouganidis:1991} for the approximation of fully nonlinear PDEs to construct difference schemes for the $\infty$-Laplacian. See also \cite{Pryer:2013} where the second author constructs an $h$-adaptive finite element scheme based on a residual error indicator and the method derived in \cite{LakkisPryer:2013}. Herein we are concerned with numerical experiments which provide further insights on the understanding of $\infty$-Harmonic mappings. To the best of our knowledge, the experiments we perform have not been attempted in the vectorial case. We consider a set of five different boundary conditions on the unit rectangle $\Om=(-1,1)^2\sub \R^2$ with target either $\R^2$ or $\R^3$ (Subsections \ref{subsection3.1}-\ref{subsection3.5}). The method we follow is based 
on the approximation of the $L^\infty$ system \eqref{1.1} by the respective $L^p$ Euler-Lagrange system, that is the $p$-Laplacian 
\beq \label{1.6}
\De_p u\, :=\, \text{div}\left(|\text{D}u|^{p-2}\text{D}u \right)\,=\,0,
\eeq
for large $p\in \N$. Our numerical scheme of choice is the finite element method. We utilise the approach described in \cite{P} for the scalar version of \eqref{1.1} where it was shown that by forming an appropriate limit we are able to select candidates for numerical approximation along a ``good'' sequence of solutions, the $p$-Harmonic mappings. 
This approach has been analytically justified by the first author in the scalar case of \eqref{1.4} in \cite{P}. Herein we justify its application to the full vectorial case of \eqref{1.1}.

The numerical method we employ is a finite element approximation, based on an earlier work of Barrett and Liu on numerical methods for elliptic systems \cite{BL}. Therein the authors prove that, for a fixed exponent $p$, the method converges to the respective $p$-Harmonic mapping under certain regularity assumptions on the solution. 
We would like to stress that significant care must be taken with numerical computations using this approach because the resulting nonlinear system is extremely badly conditioned. This owes to the nonlinearity of the problem which grows exponentially with $p$. Work to overcome this issue includes, for example, the work of Huang, Li and Liu \cite{HuangLiLiu:2007} where preconditioners based on gradient descent algorithms are designed and shown to work well for $p$ up to $1000$. We circumvent the need for such preconditioners by choosing our boundary data carefully such that $\norm{\D u} \approx 1$ over the domain. 

\emph{A word of caution:} to be perfectly clear about this paper, the purpose of our work is to demonstrate some key properties of $\infty$-Harmonic mappings by using an analytically justifiable scheme which currently is the \emph{only} technique available to give insight into the limiting vectorial problem. We note that our goal is \textbf{not} to construct an efficient approximation method for the $\infty$-Laplace system; indeed this indirect approximation of the $\infty$-Laplacian system by variational problems is not computationally efficient. 

Our results exhibit very interesting phenomena arising in the $p$-limiting case. More specifically, as $p$ increases the image of the solutions tends to ``flatten" and they behave like minimal surfaces. If the boundary condition includes two components which have ranks equal to $1,2$ (Subsections \ref{subsection3.1}, \ref{subsection3.2}), then the solutions tend to achieve the maximum possible rank throughout the domain. Moreover, as $p$ increases the angle between the $2$ partial derivative vectors appears to approach a constant value throughout without any interfaces inside the domain. If we prescribe boundary data which have rank equal to $1$ (Subsection \ref{subsection3.3}), then as $p$ increases the solutions tend to ``break" and become of rank equal to $1$ for $p=\infty$ without interfaces, while for all finite $p$ there is a region whereon the rank of the gradient is $2$ and nontrivial interfaces appear. However, in general interfaces may be formed and they may not be either smooth or with locally Euclidean topology: in subsections \ref{subsection3.4} and \ref{subsection3.5} we use as boundary data the \emph{explicit} $\infty$-Harmonic maps constructed in \cite{K2} whose interfaces are either rectangular or with triple junctions. Our results show that the $p$-Harmonic maps approach the explicit solutions as $p$ increases, forming interfaces with these shapes.

We conclude this introduction by noting that in the paper \cite{SS} the authors derived a different more singular multi-valued version of ``$\infty$-Laplacian" which describes optimal Lipschitz extensions. In our setting this amounts to changing in \eqref{1.3} from the Euclidean norm ``$|\cdot|$" we are using on $\R^{Nn}$ to the nonsmooth operator norm $\|A\|=\max_{|a|=1}|Aa|$.

\section{Basics on Vectorial Calculus of Variations in $L^\infty$ and its Fundamental Equations} \label{section2}

\subsection{$L^p$ approximations as $p\ri \infty$ of the $L^\infty$ equations} The nomenclature $\infty$-Laplacian of \eqref{1.1} owes to its very derivation as the limit of the $p$-Laplacian  \eqref{1.6} as $p\ri \infty$. In addition, the respective functionals also approximate the $L^\infty$ functional, if rescaled appropriately, in that for any $W^{1,\infty}(\Om)$ function we have
\beq \label{2.1}
E_p(u,\Om)\, :=\, \left( \int_\Om |\text{D}u|^p \right)^{1/p} \larrow \, \big\| |\text{D}u| \big\|_{L^\infty(\Om)} \, =\, E_\infty(u,\Om), \ \ \text{ as }p\ri\infty.
\eeq
In the vectorial case $N\geq 2$ the derivation of the full system \eqref{1.1} from \eqref{1.6} was first performed by the first author in \cite{K1}. We recall here the \emph{formal} derivation which is the scaffolding we employ for the numerical approximations to our solutions. Suppose $u :\Om \sub \R^n \larrow  \R^N$ is a smooth map. We rewrite the $p$-Laplacian \eqref{1.6} is index form as
\beq \label{index-form}
\sum_{i=1}^n \text{D}_i \left(|\text{D}u|^{p-2}\text{D}_i u_\al \right)\,=\,0,\ \ \ \ \al=1,...,N.
\eeq
By distributing derivatives, we have
\[
\sum_{\be=1}^N\sum_{i,j=1}^n (p-2)|\text{D}u|^{p-4} \text{D}_iu_\al\, \text{D}_ju_\be\, \text{D}^2_{ij}u_\be\, +\, |\text{D}u|^{p-2}\sum_{i=1}^n \text{D}^2_{ii}u_\al\,=\,0,\ \ \ \ \al=1,...,N
\]
and we may normalise and contract the derivatives in the first summand to find
\[
\sum_{i=1}^n  \left( \text{D}_iu_\al\, \text{D}_i\Big(\frac{1}{2}|\text{D}u|^2 \Big)\, +\, \frac{|\text{D}u|^2}{p-2} \text{D}^2_{ii}u_\al\right)\,=\,0,\ \ \ \ \al=1,...,N.
\]
In vector notation this means
\beq \label{2.2}
\text{D}u \, \text{D}\Big(\frac{1}{2}|\text{D}u|^2 \Big)\, +\, \frac{|\text{D}u|^2}{p-2}\De u\, =\, 0.
\eeq
If we let $p\ri \infty$ in \eqref{2.2}, we \emph{lose information} and we formally obtain only the system $\text{D}u\ot \text{D}u:\text{D}^2u=0$ which is one of the components of \eqref{1.1}. In the scalar case, however, this idea is correct and we obtain the full equation \eqref{1.4}. In the general vectorial case, we have the information that the two summands of \eqref{2.2} are opposite and in particular $|\text{D}u|^2\De u $ is tangential to the image $u(\Om)\sub \R^N$. In order to retain this information, for any fixed $x\in \Om$ we split $\R^N$ as the direct orthogonal sum of the range of $\text{D}u(x):\R^n \larrow  \R^N$ and of its complement
\[
\R^N\, =\, R(\text{D}u(x))\oplus R(\text{D}u(x))^\bot
\]
and by recalling \eqref{1.2}, we also set
\beq  \label{2.3}
[\text{D}u(x)]^\top\, :=\, \text{Proj}_{R(\text{D}u(x))} \, =\, I \, -\, [\text{D}u(x)]^\bot.
\eeq
By utilising \eqref{2.3} and \eqref{1.2}, we split the system \eqref{2.2} as follows
\beq \label{2.4}
\left\{ \text{D}u \, \text{D}\Big(\frac{1}{2}|\text{D}u|^2 \Big)\, +\, \frac{|\text{D}u|^2}{p-2}[\text{D}u]^\top\De u \right\}\, +\,\frac{1}{p-2}\Big\{  |\text{D}u|^2 [\text{D}u]^\bot \De u \Big\}\, =\, 0.
\eeq
Note now that the term in the first bracket is tangential since for each $x\in \Om$ it is valued in 
$R(\text{D}u(x))$, while the second term is orthogonal to the first and is valued in $R(\text{D}u(x))^\bot$. Hence, the two summands are linearly independent. We choose to renormalise \eqref{2.4} by multiplying the second summand by $p-2$. Then, after this normalisation we get
\beq \label{2.5}
  \text{D}u \, \text{D}\Big(\frac{1}{2}|\text{D}u|^2 \Big)  \, +\,  |\text{D}u|^2 [\text{D}u]^\bot \De u  \, =\, - \, \frac{|\text{D}u|^2}{p-2}[\text{D}u]^\top\De u.
\eeq
By letting $p\ri \infty$ in \eqref{2.5} we obtain the full system \eqref{1.1}. A byproduct of this derivation is that \eqref{1.1} actually can be decoupled to a pair of systems (tangential \& normal) which are independent of one another:
\beq \label{2.6}
  \text{D}u \ot \text{D}u :\text{D}^2u  \, =\, 0,\ \ \ \ \ |\text{D}u|^2[\text{D}u]^\bot\De u\, =\, 0.
\eeq

The above arguments require too much regularity and too strong estimates in order to make sense rigorously for classical $C^2$ solutions, but they are very instructive of the approach we follow. In the scalar case $N=1$ the above method of studying the $L^\infty$ equations by utilising the asymptotic limits of the $L^p$ equations as $p\ri\infty$ has proved to be very popular and successful. This idea which dates back to Aronsson has been effectively put into action by applying the theory of Viscosity Solutions to the $\infty$-Laplacian (see e.g.\ \cite{C,K} and references therein) which almost by definition is very stable under limits. Further, in view of the uniqueness theorems available in the scalar case, all subsequential limits as $p\ri \infty$ give rise to a viscosity solution of the limit equation. 

In the vectorial case $N\geq 2$ the situation is much more complicated since there is no effective counterpart of Viscosity Solutions stable under limits which would allow to prove existence with elementary estimates. Motivated partly by the equations arising in $L^\infty$, the first author has recently proposed in \cite{K7} a new duality-free theory of generalised solutions which applies to general fully nonlinear systems of any order. In particular, it allows to make sense of \eqref{1.1} in the appropriate regularity class of $W^{1,\infty}(\Om,\R^N)$ mappings. Among other results, in \cite{K7} we proved existence of solution to the Dirichlet problem for \eqref{1.1} and in \cite{K9} we proved variational characterisations of $\infty$-Harmonic maps in terms of the functional \eqref{1.3}. The idea behind this new notion of so-called \textbf{$\mD$-solutions} is briefly explained in the next subsection.

However, as it has been proved in \cite{K5} due to vectorial obstructions it is impossible to obtain uniqueness for the Dirichlet problem to \eqref{1.1} even in the class of $C^\infty$ solutions and for $n=N=2$. Evidence provided in \cite{K5} and herein suggest that \textit{the method of $L^p$ approximations as $p\ri \infty$ provides a selection principle of ``good" solution to \eqref{1.1} which has been conjectured in \cite{K7} that it is unique}.

We now state an existence result for asymptotic limits of $p$-Harmonic maps as $p\ri\infty$ needed later. The proof can be found e.g.\ in \cite{P} and is a minor vectorial extension of standard results on limits of $p$-Harmonic functions as $p\ri\infty$. Since the proof utilises only arguments involving norms, the proof in the vectorial case is essentially the same as in the scalar case given in the book \cite{K}. Let us also remind here the standard notion of weak solutions to the $p$-Laplacian: a map $u\in \sob{1}{p}_g(\W,\R^N)$ is weakly $p$-Harmonic when
\begin{equation}
  \label{eq:plap-weak}
  \int_\W \norm{\D u}^{p-2} \D u : \D v \,=\, \int_\W f \cdot v, \Foreach v \in \sob{1}{p}_0(\W,\R^N).
\end{equation}

\renewcommand{\Norm}[1]{\ensuremath{\big\|#1\big\|}}

\begin{The}[Existence of limits of $p$-Harmonic maps as $p\ri\infty$]
  \label{the:ptoinf}
  Let $\{u_p\}_{p=1}^\infty$ denote a sequence of weak
  solutions to the $p$-Laplace system \eqref{1.6} with $u_p\in \sob{1}{p}_g(\W,\R^N)$. Then, there exists a subsequence such  that as $p\ri \infty$ that sequence converges uniformly to a mapping $u_\infty\in \sob{1}{\infty}(\W,\R^N)$. Namely,
  \begin{equation}
    u_{p_j} \larrow  u_\infty \text{ in } C^{0}(\overline{\Omega},\R^N), \ \ \text{ as }p\ri\infty.
  \end{equation}
\end{The}
We emphasise that the map $u_\infty$ is a \textbf{candidate} $\infty$-Harmonic mapping, that is a generalised solution to \eqref{1.1}. This is true in the scalar case $N=1$ in the sense of Viscosity Solutions of Crandall-Ishii-Lions. In the vectorial case, it is true as well at least in the case of the ODE system arising from \eqref{1.5} when $n=1$ in the sense of $\mD$-solutions of the first author (see \cite{K8}). We conjecture to be true in the case of \eqref{1.1} as well, but this is not a consequence of the current results of \cite{K7} since the method of the existence proof was based on an ad-hoc method (an analytic counterpart of Gromov's ``Convex Integration" for a differential inclusion) rather than on $p$-Harmonic approximations. A complete proof of this conjecture, at least to date, eludes us, but recently we have made significant progress in this regard.

\subsection{Generalised solutions of the $L^\infty$ equations} For the sake of completeness we briefly motivate here the definition of generalised solutions to \eqref{1.1}. Since we do not utilise it in an essential manner in this paper, we refrain from giving all the details which can be found in \cite{K7} and \cite{K8}-\cite{K11}. The idea applies to general fully nonlinear systems of any order and allows for merely measurable solutions. It is based on the following observation: a map $u : \Om\sub \R^n \larrow  \R^N$ in $C^2(\Om,\R^N)$ is a solution to \eqref{1.1} if and only if for any compactly supported function $\Phi \in C^0_c(\R^{Nn^2}_s)$, we have
\beq \label{2.7}
\int_{ \R^{Nn^2}_s }\Phi(\X) \, \Big(\text{D}u\ot \text{D}u +|\text{D}u|^2[\text{D}u]^\bot \!\ot I\Big) :\X\, d [\de_{\text{D}^2u} ](\X)\,=\,0,\ \ \ \text{on } \Om.
\eeq 
The equation \eqref{2.7} is a mere restatement of the system \eqref{1.1} where we just change the viewpoint and instead of considering the hessian as a classical map $\text{D}^2u :\Om\sub \R^n \larrow  \R^{Nn^2}_s$, we instead view it as a probability valued map given by the Dirac mass at the hessian:
\[
x\lmapsto \de_{\text{D}^2u(x)}\ : \ \Om\sub \R^n \larrow  \mathscr{P}(\R^{Nn^2}_s).
\]
This change of viewpoint turns out to be very fruitful since by attaching one point and compactifying $\R^{Nn^2}_s$, we have that for any perhaps nonsmooth map, there always exist limits of the difference quotients in the appropriate space of probability-valued maps which may not be the concentration measures $\de_{\text{D}^2u}$ but instead more general probability valued maps $\mD^2u : \Om\sub \R^n \larrow  \mathscr{P}(\R^{Nn^2}_s\cup \{\infty\})$ called \textbf{diffuse hessians}. More precisely, if $D^{1,h}$ denotes the first difference quotient operator on $\R^n$, our generalised hessians are the subsequential limits of the form
\[
\de_{D^{1,h}\text{D}u}\, \weakstar\, \mD^2u, \ \text{ as $h\ri 0$ in the Young measures }\ \Om\sub \R^n\larrow \R^{Nn^2}_s\cup \{\infty\}.
\]
The respective notion which generalises \eqref{2.7} is called \textbf{$\mD$-solutions} to the $\infty$-Laplacian \eqref{1.1} and is the primary notion of generalised solution in this context for the vectorial case.

\subsection{Some explicit smooth $\infty$-Harmonic mapping in $2\by2$ dimensions} 
\begin{figure}[!h]
  \caption[Numerical Results for Problem \eqref{eqn:Problem:1} with
  $P^1$ elements]
  {\label{fig:1}
     
  }
  \begin{center}
    \subfloat[{\label{fig:1c}
        {The graph of the parametrisation function $K$ of (I) used in \ref{2.8}.}
    }]{
      \includegraphics[scale=\figscale,width=0.25\figwidth]{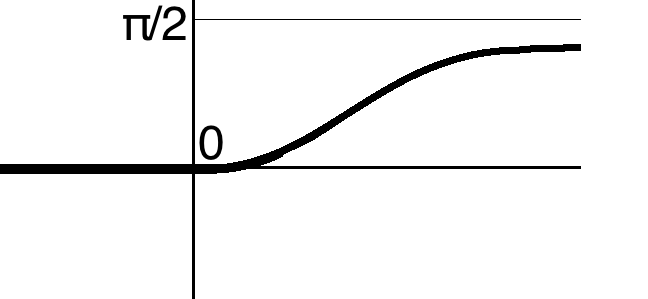}
    }
    \qquad 
    \qquad 
    \qquad 
    \qquad 
    \qquad 
    \qquad 
    \subfloat[{\label{fig:1d}
        {The graph of the parametrisation function $K$ of (II) used in \ref{2.8}.}
    }]{
      \includegraphics[scale=\figscale,width=0.3\figwidth]{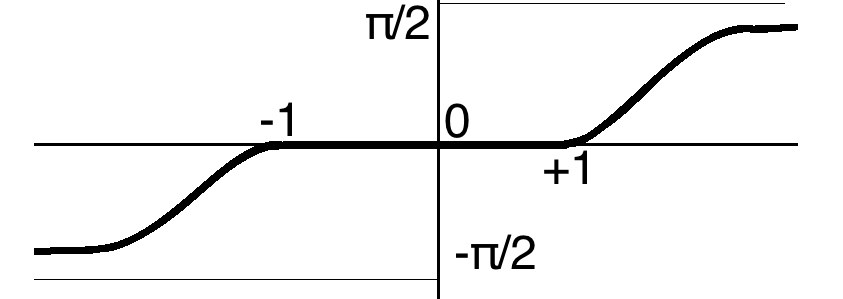}
    }
    \\
    \subfloat[{
        {Phases and interfaces of the smooth $\infty$-Harmonic map given by the explicit formula \ref{2.8} when $K$ is as in (I).}
    }]{\label{fig:1a}
      \includegraphics[scale=\figscale,width=0.25\figwidth]{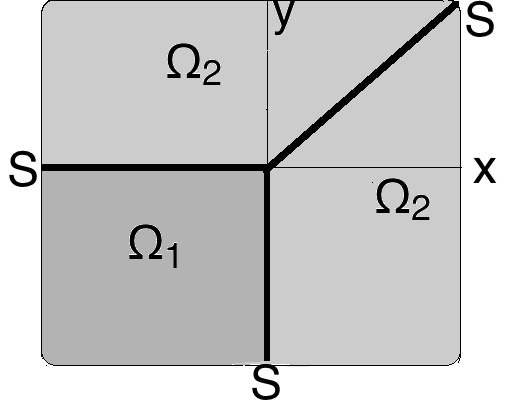}
    }
    \qquad 
    \qquad 
    \qquad 
    \qquad 
    \qquad 
    \qquad 
    \subfloat[{\label{fig:1b}
        {Phases and interfaces of the smooth $\infty$-Harmonic map given by the explicit formula \ref{2.8} when $K$ is as in (II).}
    }]{
      \includegraphics[scale=\figscale,width=0.25\figwidth]{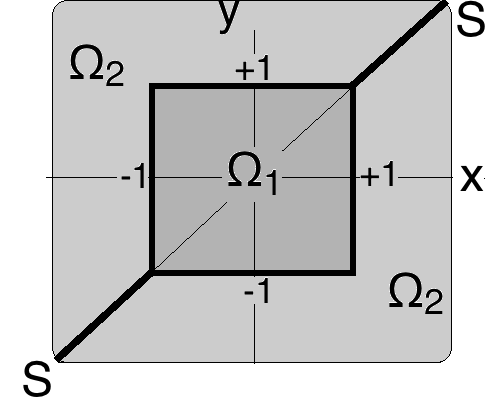}
    }
  \end{center}
\end{figure}

The following explicit solutions of the system \eqref{1.1} have been constructed in \cite{K2} and we briefly recall them here because they are utilised later in Section \ref{section3}. Let $u : \R^2 \larrow  \R^2$ be a map in $C^1(\R^2,\R^2)$. We set
\beq \label{2.11}  
\left\{\ \
\begin{split}
\Om_2 \, &:=\, \Big\{x\in \R^2\,:\,\rk\big(\text{D}u(x)\big)=2\Big\},\\
 \Om_1\, &:=\,  \inter  \Big\{x\in \R^2\,:\, \rk\big(\text{D}u(x)\big)\leq 1\Big\},\\
 \S\, & :=\, \R^2 \set (\Om_1 \cup \Om_2).
\end{split}
\right.
\eeq
Here ``rk" denotes the rank and ``int" denotes the topological interior. We call \emph{$\Om_2$ the $2$-dimensional phase, $\Om_1$ the $1$-dimensional phase and $\S$ the interface} of the map $u$. Obviously, $\R^2=\Om_2 \cup \Om_1 \cup \S$. On $\Om_2$ $u$ is a local diffeomorphism and on $\Om_1$ it is ``essentially scalar". In \cite{K2} we proved that the explicit formula
\beq  \label{2.8}
u(x,y)\ :=\ \int_y^x e^{iK(t)}\, \text{d}t
\eeq
defines a smooth explicit $\infty$-Harmonic map $(-1,1)^2\sub\R^2 \larrow \R^2$ when $K\in C^1(\R)$ and $\|K\|_{C^0(\R)} <\pi/2$. Obviously, $e^{it}$ symbolises $(\cos t, \sin t)^\top$. Moreover: 

\noi (I) If $K$ is qualitatively as in Figure \ref{fig:1c}, namely $K\equiv 0$ on $(-\infty, 0]$ and $K'>0$ on $(0,\infty)$, then $u$ is affine on $\Om_1$ and $\Om_1$, $\Om_2$, $\S$ are as in Figure \ref{fig:1a}, i.e.
\[
\Om_1 = \{x,y < 0\},\ \ \S = \p \Om_1 \cup \{x=y\geq 0\}, \ \ \Om_2 = \R^2 \set (\Om_1 \cup \S).
\]
\noi (II)  If $K$ is qualitatively as in Figure \ref{fig:1d}, namely if $K\equiv 0$ on $[-1,+1]$ and $K'>0$ on $(-\infty, -1) \cup (1,\infty)$, then $u$ is affine on $\Om_1$ and $\Om_1$, $\Om_2$, $\S$ are as in Figure \ref{fig:1b}, i.e.
\[
\Om_1 = \{-1< x,y <1 \},\ \ \S = \p \Om_1 \cup \{x=y, |y|\geq 1\}, \ \ \Om_2 = \R^2 \set (\Om_1 \cup \S).
\]

The above examples show that in the vectorial case very complicated phenomena can arise even for smooth solutions. We will further examine these phenomena in Section \ref{section4} numerically by studying $p$-Harmonic mappings for increasing values of $p$ using boundary data provided by \ref{2.8} with choices of $K$ as in (I), (II).

\section{Numerical Approximations of $\infty$-Harmonic mappings for $n=N=2$ and $n=2 < N=3$} \label{section3}

In this section we describe the technique we use to approximate $\infty$-Harmonic mappings. The method we use is a conforming finite element discretisation of the $p$-Laplacian analysed in \cite{BL} for fixed $p$. We will describe the discretisation and justify its application to the problem at hand by studying the behaviour as the meshsize parameter tends to zero and as $p$ gets large. To that end we will also extend the results given in \cite{P} to the vectorial case.

We let $\T{}$ be an admissible triangulation of $\W$,
namely, $\T{}$ is a finite collection of sets such that

\smallskip
\begin{enumerate}
\item $K\in\T{}$ implies $K$ is an open triangle ,
\item for any $K,J\in\T{}$ we have that $\closure K\meet\closure J$ is
  either $\emptyset$, a vertex, an
  edge, or the whole of $\closure K$ and $\closure J$ and
\item $\union{K\in\T{}}\closure K=\closure\W$.
\end{enumerate}

\smallskip
\noi The shape regularity constant of $\T{}$ is defined as the number
\begin{equation}
  \label{eqn:def:shape-regularity}
  \mu(\T{}) \,:=\, \inf_{K\in\T{}} \frac{\rho_K}{h_K},
\end{equation}
where $\rho_K$ is the radius of the largest ball contained inside
$K$ and $h_K$ is the diameter of $K$. An indexed family of
triangulations $\setof{\T n}_n$ is called \emph{shape regular} if 
\begin{equation}
  \label{eqn:def:family-shape-regularity}
  \mu\, :=\, \inf_n\mu(\T n)\, >\, 0.
\end{equation}
Further, we define $\funk h\W \R$ to be the {piecewise
  constant} \emph{meshsize function} of $\T{}$ given by
\begin{equation}
h\, \equiv \, h(\vec{x})\, :=\, \max_{\closure K\ni \vec{x}}\, h_K.
\end{equation}
A mesh is called quasiuniform when there exists a positive constant
$C$ such that $\max_{x\in\Omega} h \le C \min_{x\in\Omega} h$. In what
follows we shall assume that all triangulations are shape-regular and
quasiuniform.


We let $\poly k(\T{})$ denote the space of piecewise polynomials of
degree $k$ over the triangulation $\T{}$,\ie
\begin{equation}
  \poly k (\T{}) \, = \, \big\{ \phi \text{ such that } \phi|_K \in \poly k (K) \big\}
\end{equation}
 and introduce the \emph{finite element space}
\begin{gather}
  \label{eqn:def:finite-element-space}
  \fes \,:= \, \poly k(\T{}) \cap \cont{0}(\W)
\end{gather}
to be the usual space of continuous piecewise polynomial
functions of degree $k$.



\begin{Defn}[$\leb{2}(\W)$ projection operator]
  \label{eq:l2proj}
  The $\leb{2}(\W)$ projection operator, $P_h : \leb{2}(\W) \larrow  \fes$ is
  defined for $v\in\leb{2}(\W)$ such that
  \begin{equation}
    \int_\W P_h v \, \Phi \,=\, \int_\W v \, \Phi,  \Foreach \phi\in\fes.
  \end{equation}
  It is well known that this operator satisfies the following
  approximation properties for $v \in \sob{1}{p}(\W)$ 
  \begin{gather}
    \lim_{h\ri  0}\Norm{v - P_h v}_{\leb{p}(\W)} =\, 0
    \\
    \lim_{h\ri  0}\Norm{\D v - \D \qp{ P_h v}}_{\leb{p}(\W)} = \, 0.
  \end{gather}
\end{Defn}

\subsection{Galerkin discretisation}
We consider the Galerkin discretisation of (\ref{index-form}), to find $U \in \qb{\fes}^N$ with $U\vert_{\partial \W} = P_h g$ such that
\begin{equation}
  \label{eq:plapdis}
  \int_\W \norm{ \D U}^{p-2} \D U : {\D \Phi} \,=\, 0,\ \ \Foreach \Phi \in \qb{\fes}^N.
\end{equation}
This is a conforming finite element discretisation of the vectorial $p$-Laplacian system proposed in \cite{BL}. 

\begin{Pro}[existence and uniqueness of solution to (\ref{eq:plapdis})]
  There exists a unique solution of both the weak formulation (\ref{eq:plap-weak}) and the Galerkin approximation (\ref{eq:plapdis}).
\end{Pro}
\begin{Proof}.
  Existence and uniqueness of this problem follows from examination of the $p$-functional
  \begin{equation}
    \label{eq:p-energy}
    E_p(u,\W) \, =\, \left(\int_\Om |\text{D}u|^p \right)^{1/p}. 
  \end{equation}
  Notice that (\ref{eq:p-energy}) is strictly convex and coercive on $\sob{1}{p}_0(\W,\R^N)$ so we may apply standard arguments from the Calculus of Variations showing that the minimisation problem is well posed. Hence, there exists $u\in \sob{1}{p}_g(\W,\R^N)$ such that
  \begin{equation}
    \label{eq:min}
    E_p(u,\W) \,=\, \min_{ v\in \sob{1}{p}_0(\W,\R^N) } E_p(v,\W).
  \end{equation}
Noticing that the weak problem (\ref{eq:plap-weak}) is the weak Euler-Lagrange equation for (\ref{eq:min}) we have equivalence of (\ref{eq:plap-weak}) and (\ref{eq:min}) as such the weak formulation is also well posed. To see this for the Galerkin approximation notice that $\qb{\fes}^N \subset  \sob{1}{p}(\W,\R^N)$. As such the minimisation problem over $\qb{\fes}^N$ is equivalent to (\ref{eq:plapdis}) and the same argument applies as in the continuous case.
\end{Proof}

\begin{The}[convergence of the discrete scheme to weak solutions]
  \label{the:convergence-of-scheme}
  For fixed $p$ let $\{U_p\}$ be the finite element approximation generated by
  solving (\ref{eq:plapdis}) (indexed by $h$) and $u_p$, the weak solution of (\ref{eq:plap-weak}), then we have that
  \begin{equation}
    U_p \larrow  u_p  \ \text{ in }  \cont{0}(\overline{\W},\R^N), \text{ as } h \ri  0.
  \end{equation}
\end{The}
\begin{Proof}
  We begin by noting the discrete weak formulation (\ref{eq:plapdis})
  is equivalent to the minimisation problem: Find $U\in\qb{\fes}^N$ with $U|_{\partial \W} = P_h g$ such that
  \begin{equation}
    \label{eq:dismin}
    E_p(U, \W) \,=\, \min_{V\in\qb{\fes}^N} E_p(V,\W).
  \end{equation}
  Using this, we immediately have
  \begin{equation}
    \Norm{\norm{\D U}}_{\leb{p}(\W)}^p \,\leq \, E_p(U, \W) \,\leq \, E_p(P_h g, \W) \,\leq\, \Norm{\norm{\D \qp{P_h g}}}_{{\leb{p}(\W)}}^p.
  \end{equation}
  In view of the stability of the $\leb{2}$ projection in
  $\sob{1}{p}(\W)$ \cite{CrouzeixThomee:1987} we have
  \begin{equation}
    \Norm{\norm{\D U}}_{{\leb{p}(\W)}} \,\leq \, C,
  \end{equation}
  uniformly in $h$.
  Hence by weak compactness there exists a (weak) limit to the finite
  element sequence, which we will call $u^*$. \text{D}ue to the weak
  semicontinuity of $E_p(\cdot, \W)$ we have 
  \begin{equation}
    E_p(u^*, \W) \, \leq \, E_p(U, \W).
  \end{equation}
  In addition, in view of the approximation properties of $P_h$ given in
  Definition \ref{eq:l2proj} we have for any $v\in \cont{\infty}(\W,\R^N)$ that
  \begin{equation}
    E_p(v, \W) \,= \, \liminf_{h\ri  0} E_p(P_k v, \W).
  \end{equation}
  Using the fact that $U$ is a discrete minimiser of (\ref{eq:dismin})
  we have
  \begin{equation}
    E_p(u^*, \W) \, \leq \, E_p(U, \W) \, \leq \, E_p(P_h v, \W),
  \end{equation}
  whence sending $h\larrow  0$ we see 
  \begin{equation}
    E_p(u^*, \W) \, \leq \, E_p(v, \W).
  \end{equation}
  Now, as $v$ was generic we may use density arguments and that $u_p$ was the unique minimser to conclude $u^* = u_p$, concluding the proof.
\end{Proof}

\begin{The}[convergence]
  \label{the:convergence}  
  Let $U_p$ be the Galerkin solution of (\ref{eq:plapdis}) and $u_\infty$ the candidate $\infty$-Harmonic mapping. Then, along a subsequence we have
  \begin{equation}
    U_{p_j} \larrow  u_\infty \text{ in } \cont{0}, \text{ as } p \ri  \infty \text{ and } h \ri  0.
  \end{equation}
\end{The}
\begin{Proof}
  The proof consists of combining Theorems \ref{the:ptoinf} and \ref{the:convergence-of-scheme} and noticing that
  \begin{equation} 
    \Norm{\norm{U_{p_j} - u_\infty}}_{\cont{0}(\W)}
   \, \leq \,
    \Norm{\norm{U_{p_j} - u_{p_j}}}_{\cont{0}(\W)}
   \, +
\,    \Norm{\norm{u_{p_j} - u_\infty}}_{\cont{0}(\W)}.
  \end{equation}
\end{Proof}


\section{Numerical experiments}
\label{section4}

In this section we summarise extensive numerical experiments which focus on quantifying the structure of solutions to the $\infty$-Laplacian PDE system \eqref{1.1}. This is achieved using Galerkin approximations to the $p$-Laplacian for sufficiently high values of $p$. We focus on studying the behaviour solutions have as $p$ increases which allow us to make various conjectures on the behaviour of their asymptotic limit as $p\ri\infty$.

\begin{Rem}[practical computation of (\ref{eq:plapdis}) for large $p$]
  The computation of $p$-Harmonic mappings is an extremely challenging problem in its own right. The class of nonlinearity in the problem results in the algebraic system, which ultimately yields the finite element solution, being extremely badly conditioned. One method to tackle this class of problems is by using preconditioners based on descent algorithms \cite{HuangLiLiu:2007}. For extremely large $p$, say $p \geq 10000$, this may be required; however for our purposes we restrict our attention to $p \sim 1000$. This yields sufficient accuracy for the results we want to illustrate.  
  
We emphasise that even the case $p\sim 1000$ is computationally tough to handle. The numerical approximation we are using is based on a Newton solver. As it is well known, Newton solvers require a sufficiently close initial guess in order to converge. A reasonable initial guess for the $p$-Laplacian is given by numerically approximating with the $q$-Laplacian for $q < p$ sufficiently close to $p$. This leads to an iterative process in the generation of the initial guess, i.e., we solve the $2$-Laplacian as an initial guess to the $3$-Laplacian which serves as an initial guess to the $4$-Laplacian, and so on. 
\end{Rem}

To tie into the explicit examples given in Figure \ref{fig:1} we are particularly interested in the rank of the solution. Except for the intrinsic interest, this relates directly to a deeper understanding of the solutions to $\infty$-Laplace system since the coefficients of \eqref{1.1} are discontinuous on the interfaces of the solution. We compute this by calculating $\det{\D U}$ and representing the areas $\Om_2$, $\Om_1$ and $\S$ of \eqref{2.11} by plotting contours of the function $\det{\D U}$.

\subsection{Solutions $(-1,+1)^2\sub \R^2\larrow  \R^2$ with mixed rank-two and rank-one boundary data.} \label{subsection3.1}

\begin{figure}[!h]
  \caption[]
  {\label{fig:mixed-r1-r2}
    An illustration of the rank of the solution to the vectorial $p$-Laplacian with the mixed rank-one and rank-two boundary conditions given in Section \ref{subsection3.1} for various increasing values of $p$. Here we plot $\det{\D U}$ and associated contour lines. These are plotted at increments of $0.05$. Notice as $p$ increases, the region where the solution is not of full rank, $\W_1$, decreases in size.
  }
  \begin{center}
    \subfloat[{\label{fig:r1r2a1}
        $\det{\D U}$ for the vectorial $2$-Laplacian.
    }]{
      \includegraphics[scale=\figscale,width=0.35\figwidth]{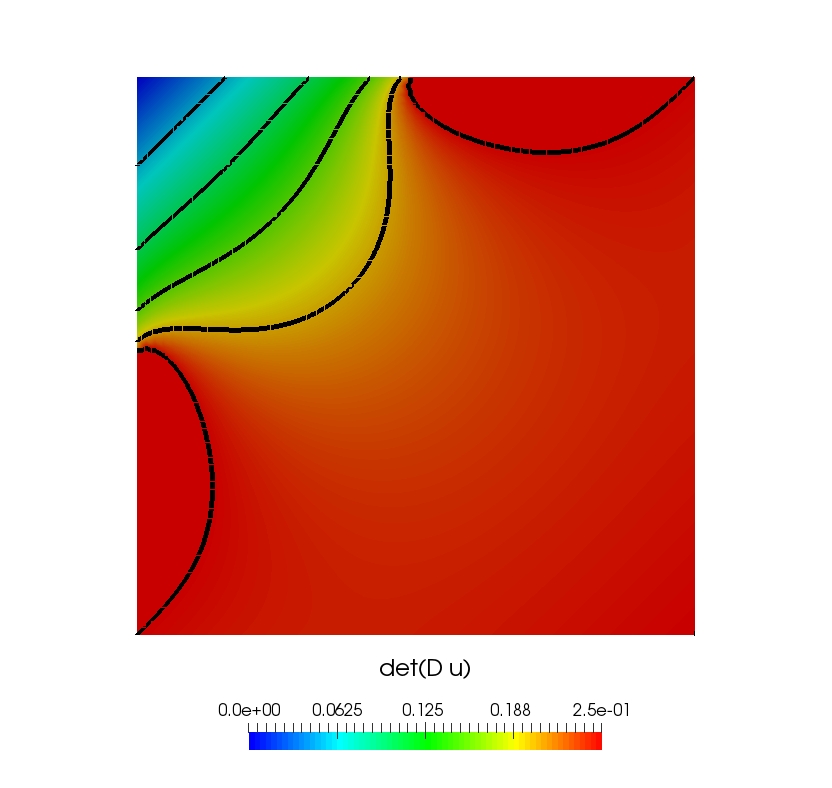}
    }
    \hfill
    \subfloat[{\label{fig:r1r2a2}
        $\det{\D U}$ for the vectorial $8$-Laplacian.
    }]{
      \includegraphics[scale=\figscale,width=0.35\figwidth]{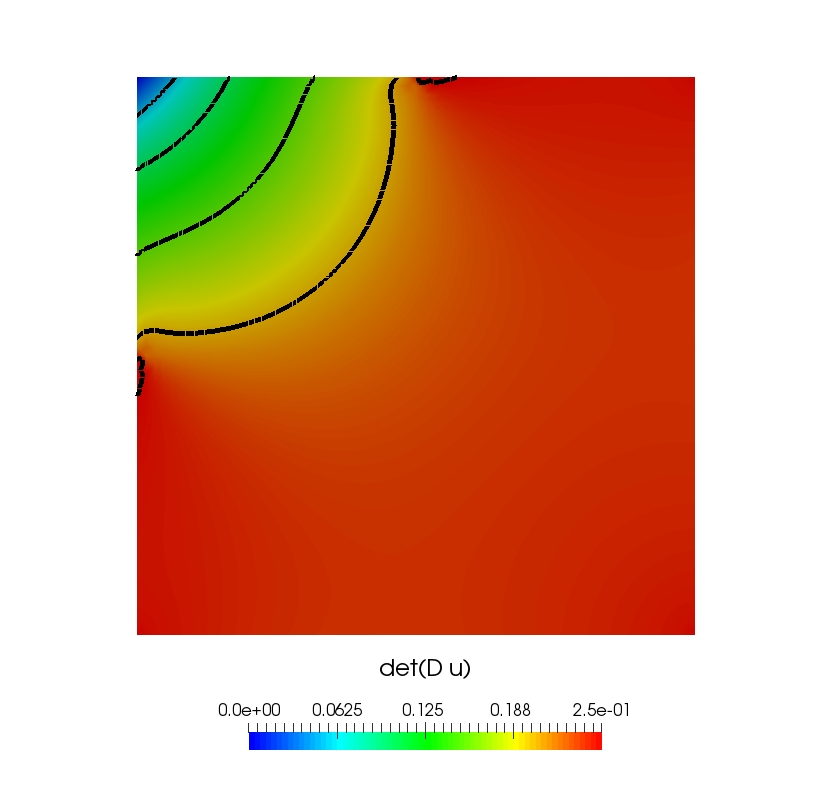}
    }
    \hfill
    \subfloat[{\label{fig:r1r2a3}
        $\det{\D U}$ for the vectorial $64$-Laplacian.
    }]{
      \includegraphics[scale=\figscale,width=0.35\figwidth]{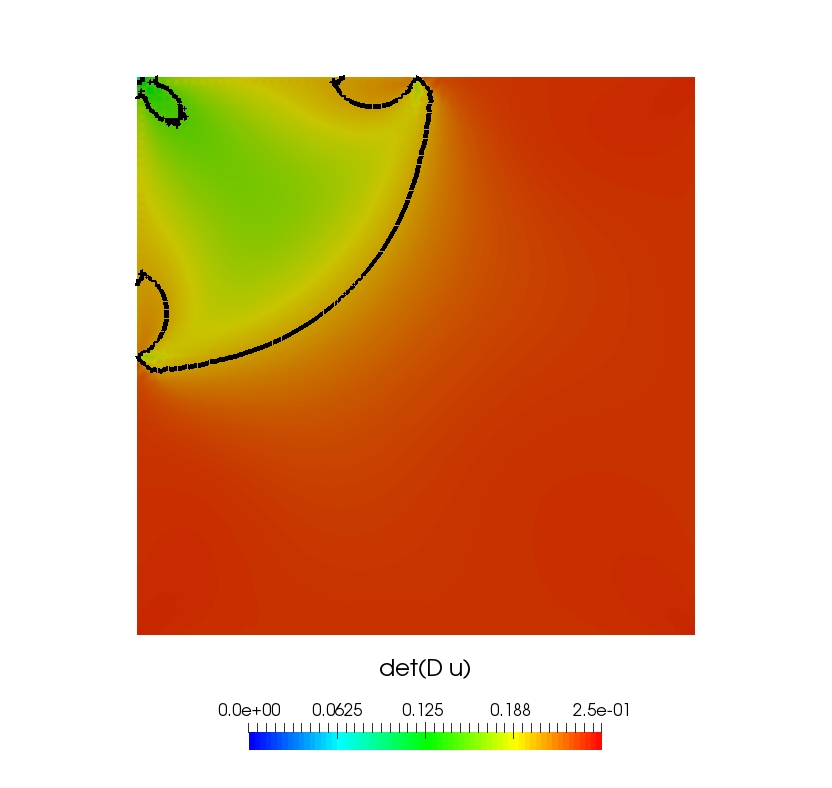}
    }
    \hfill
    \subfloat[{\label{fig:r1r2a4}
        $\det{\D U}$ for the vectorial $256$-Laplacian.
    }]{
      \includegraphics[scale=\figscale,width=0.35\figwidth]{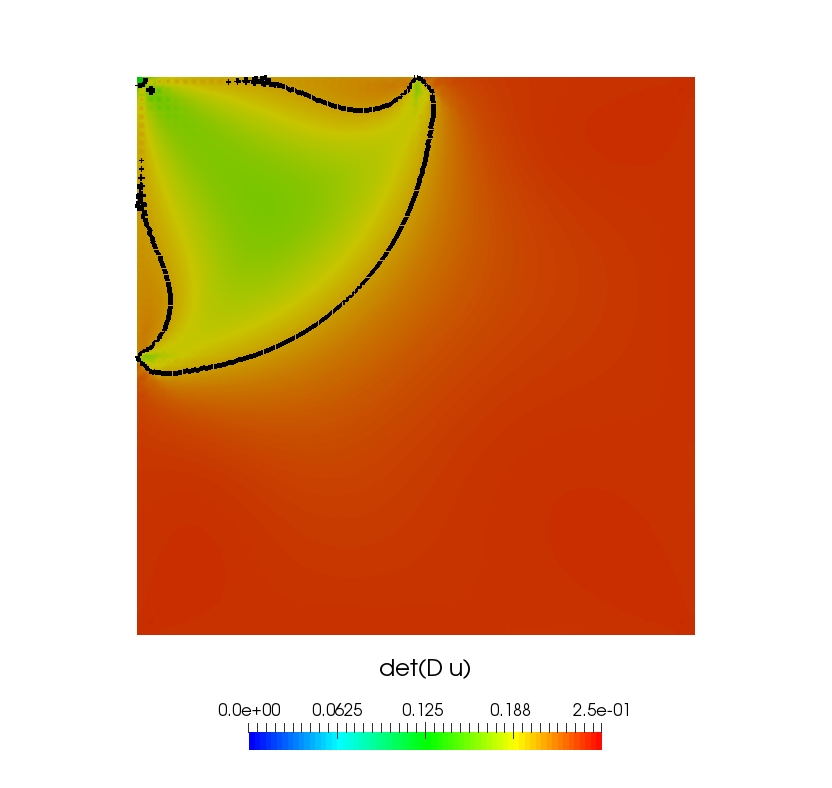}
    }
    \hfill
    \subfloat[{\label{fig:r1r2a5}
        $\det{\D U}$ for the vectorial $512$-Laplacian.
    }]{
      \includegraphics[scale=\figscale,width=0.35\figwidth]{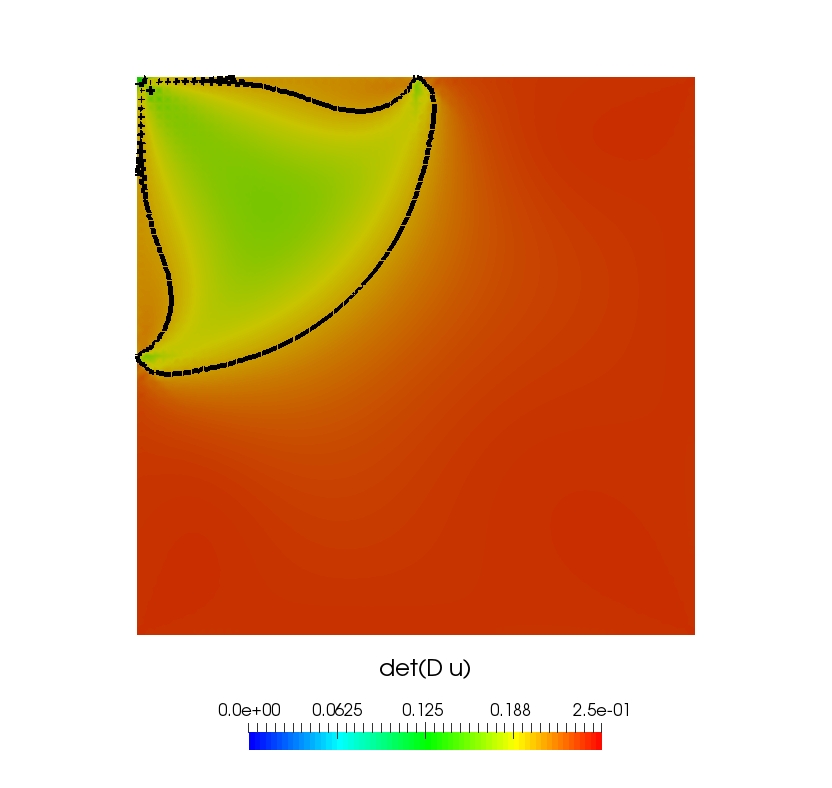}
    }
    \hfill
    \subfloat[{\label{fig:r1r2a6}
        $\det{\D U}$ for the vectorial $1024$-Laplacian.
    }]{
      \includegraphics[scale=\figscale,width=0.35\figwidth]{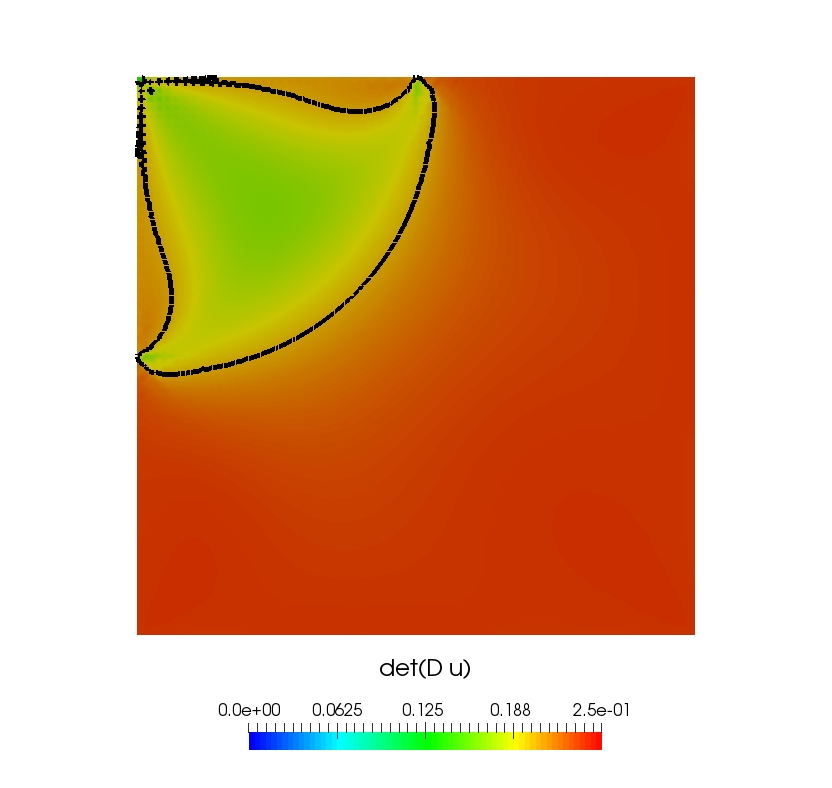}
    }
    \end{center}
  \end{figure}

In this test we construct approximations of solutions of \eqref{1.1} with mixed rank-two and rank-one boundary data. We take
\begin{equation}
  g(x,y)\, :=
\,
\left\{
  \begin{array}{ll}
    \dfrac{1}{2} \big(x,y\big)^\top , & \text{ if } x \geq 0 \text{ or } y \leq 0, \ms
    \\
    \dfrac{1}{4}\big(x+y-1, x+y+1\big)^\top,  & \text{ otherwise.}
  \end{array}
  \right.
\end{equation}
This gives us rank-one data in the quadrant $x<0$ and $y<0$ and rank-two data elsewhere. The results are illustrated in Figure \ref{fig:mixed-r1-r2}.

\subsection{Solutions $(-1,+1)^2\sub \R^2\larrow  \R^3$ with mixed rank-two and rank-one boundary data.} \label{subsection3.2}

In this experiment we examine an extension to the example given in Section \ref{subsection3.1} to the case $N=3$. We study the image of the problem with boundary data given as
\begin{equation}
  g(x,y) \,:=
  \,
\left\{
  \begin{array}{ll}
    \dfrac{1}{2} \big(x,y,x\big)^\top, & \text{ if } x \geq 0 \text{ or } y \leq 0, \ms
    \\
    \dfrac{1}{4}\big(x+y-1, x+y+1, x+y-1\big)^\top , &\text{ otherwise.}    
  \end{array}
  \right.
\end{equation}

\begin{figure}[!h]
  \caption[]
  {\label{fig:3d}
    The image of the solution to the vectorial ($N=3$) $p$-Laplacian with mixed rank-one and rank-two boundary conditions given in Section \ref{subsection3.2} for two values of $p$. 
  }
  \begin{center}
    \subfloat[{\label{fig:3da1}
        The image of the vectorial $2$-Laplacian.
    }]{
      \includegraphics[scale=\figscale,width=0.4\figwidth]{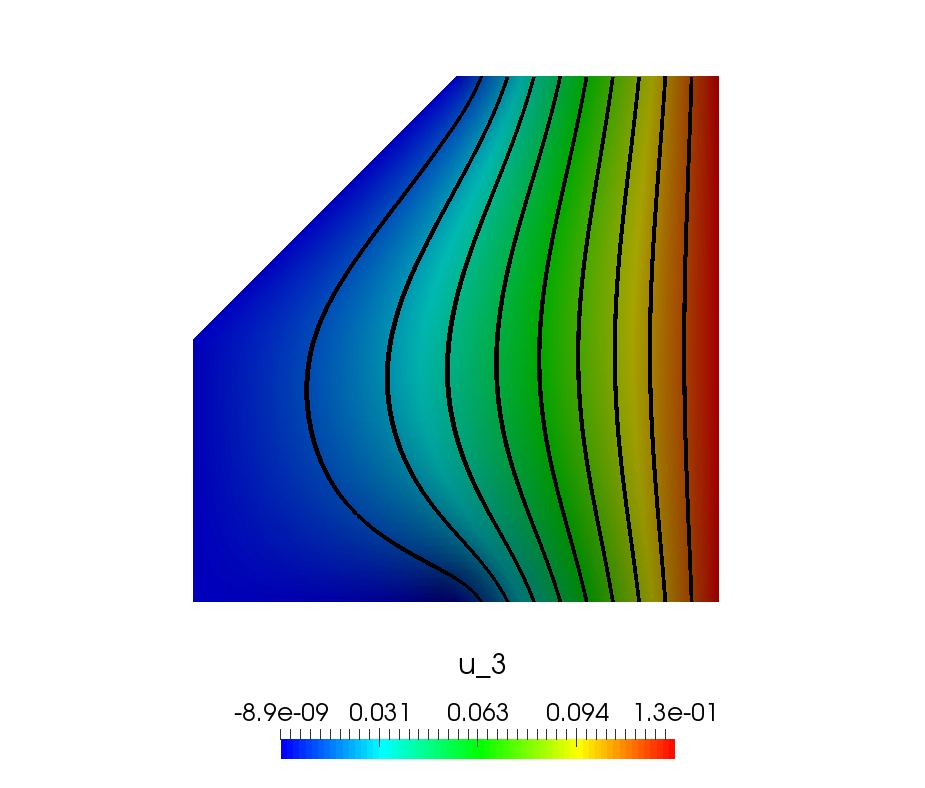}
    }
    \hfill
    \subfloat[{\label{fig:3da3}
        The image of the vectorial $2$-Laplacian, side view.
    }]{
      \includegraphics[scale=\figscale,width=0.4\figwidth]{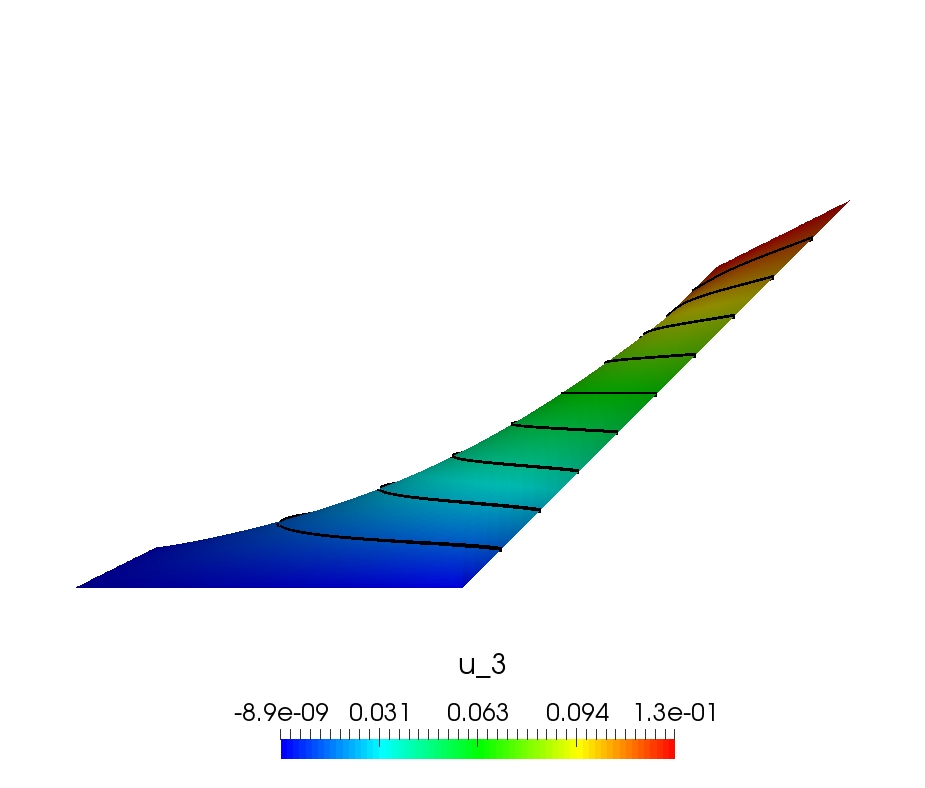}
    }
    \hfill
    \subfloat[{\label{fig:3da7}
        The image of the vectorial $820$-Laplacian.
    }]{
      \includegraphics[scale=\figscale,width=0.4\figwidth]{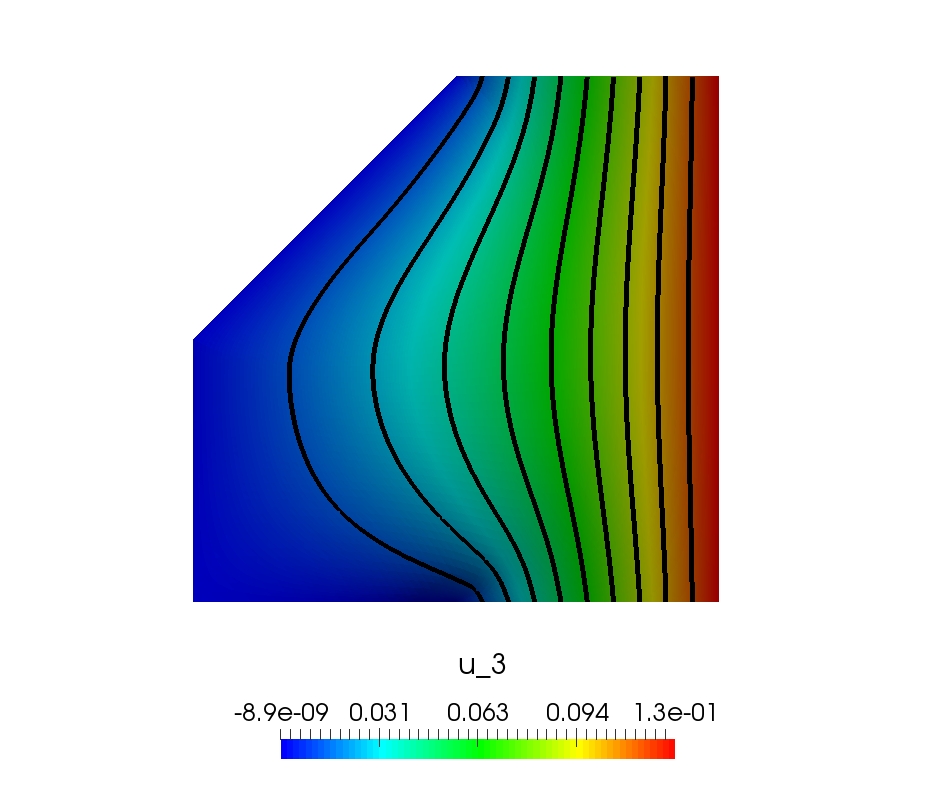}
    }
    \hfill
    \subfloat[{\label{fig:3da7}
        The image of the vectorial $820$-Laplacian.
    }]{
      \includegraphics[scale=\figscale,width=0.4\figwidth]{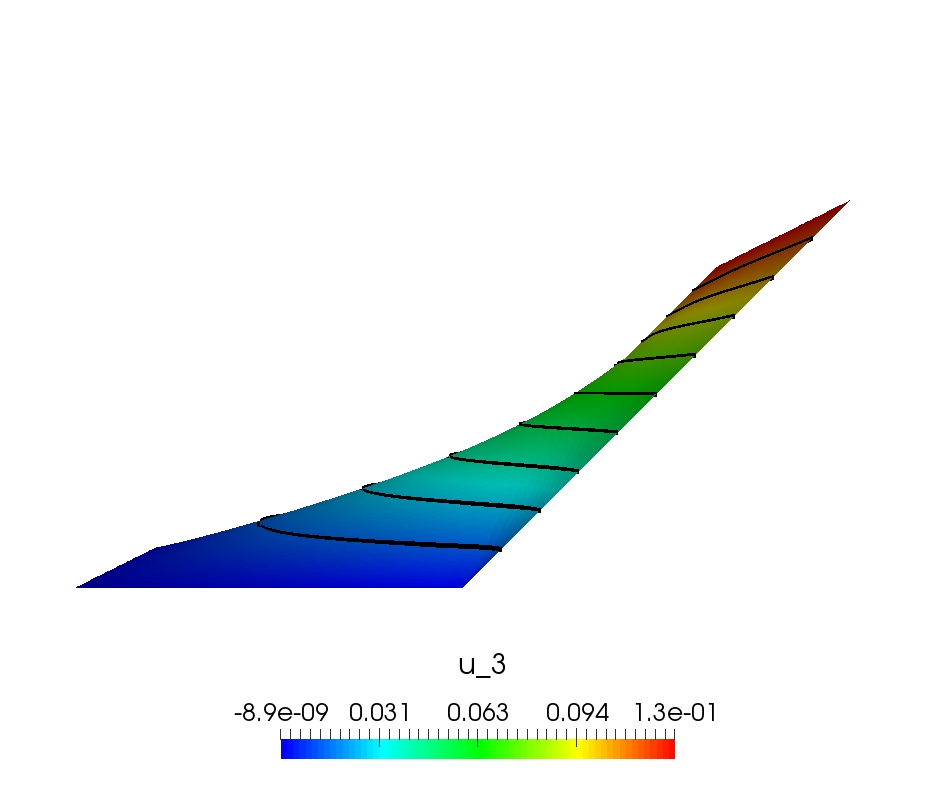}
    }
    \hfill
    \subfloat[{\label{fig:3da7}
        A comparison of the contours. The $2$-Laplacian is the lighter colour (blue).
    }]{
      \includegraphics[scale=\figscale,width=0.4\figwidth]{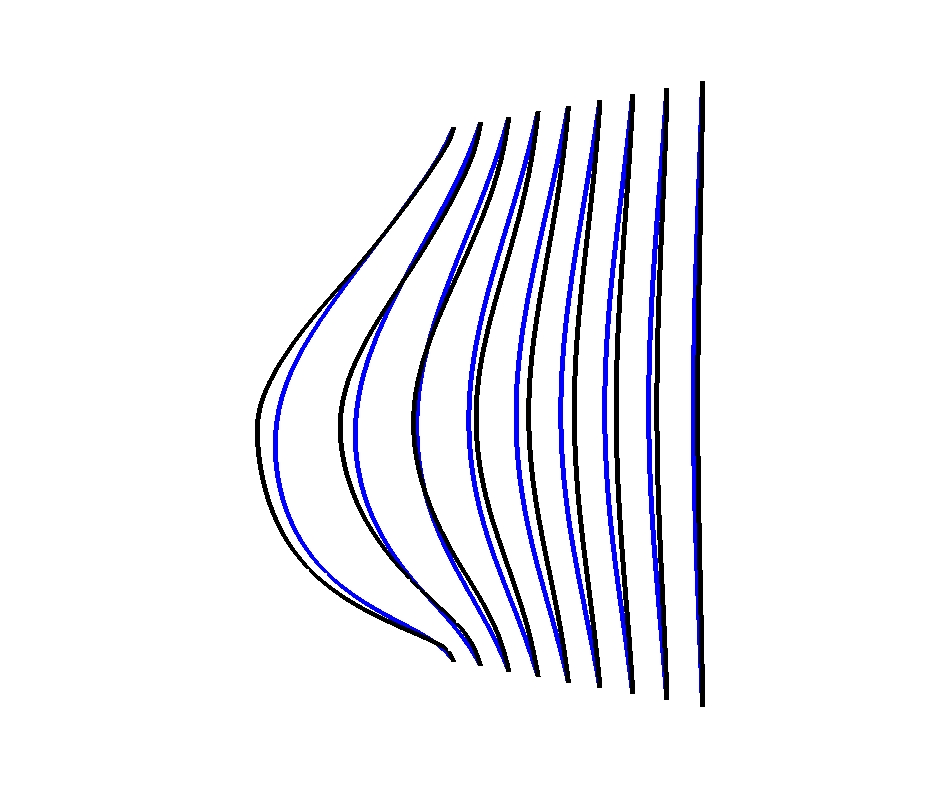}
    }
    \hfill
    \subfloat[{\label{fig:3da7}
        A comparison of the contours from the side view. The $2$-Laplacian is the lighter colour (blue).
    }]{
      \includegraphics[scale=\figscale,width=0.4\figwidth]{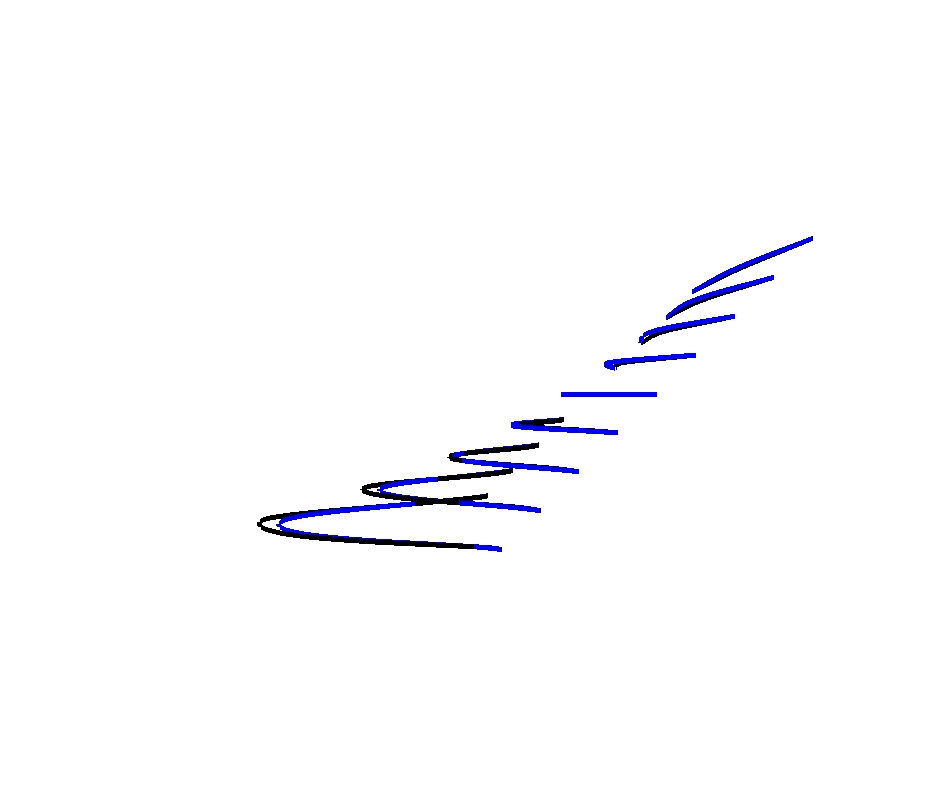}
    }
  \end{center}
\end{figure}

\subsection{Solutions $(-1,+1)^2\sub \R^2\larrow  \R^2$ with rank-one boundary data.} \label{subsection3.3}

\begin{figure}[!h]
  \caption[]
  {\label{fig:rank1im}
    The image of the solution to the vectorial $p$-Laplacian with the rank-one boundary conditions given in (\ref{eq:rank1-bcs}) for various increasing values of $p$. Notice as $p$ increases, the image flattens. It is behaving like a minimal surface.
  }
  \begin{center}
    \subfloat[{\label{fig:r1ia1}
        The image of the vectorial $2$-Laplacian.
    }]{
      \includegraphics[scale=\figscale,width=0.35\figwidth]{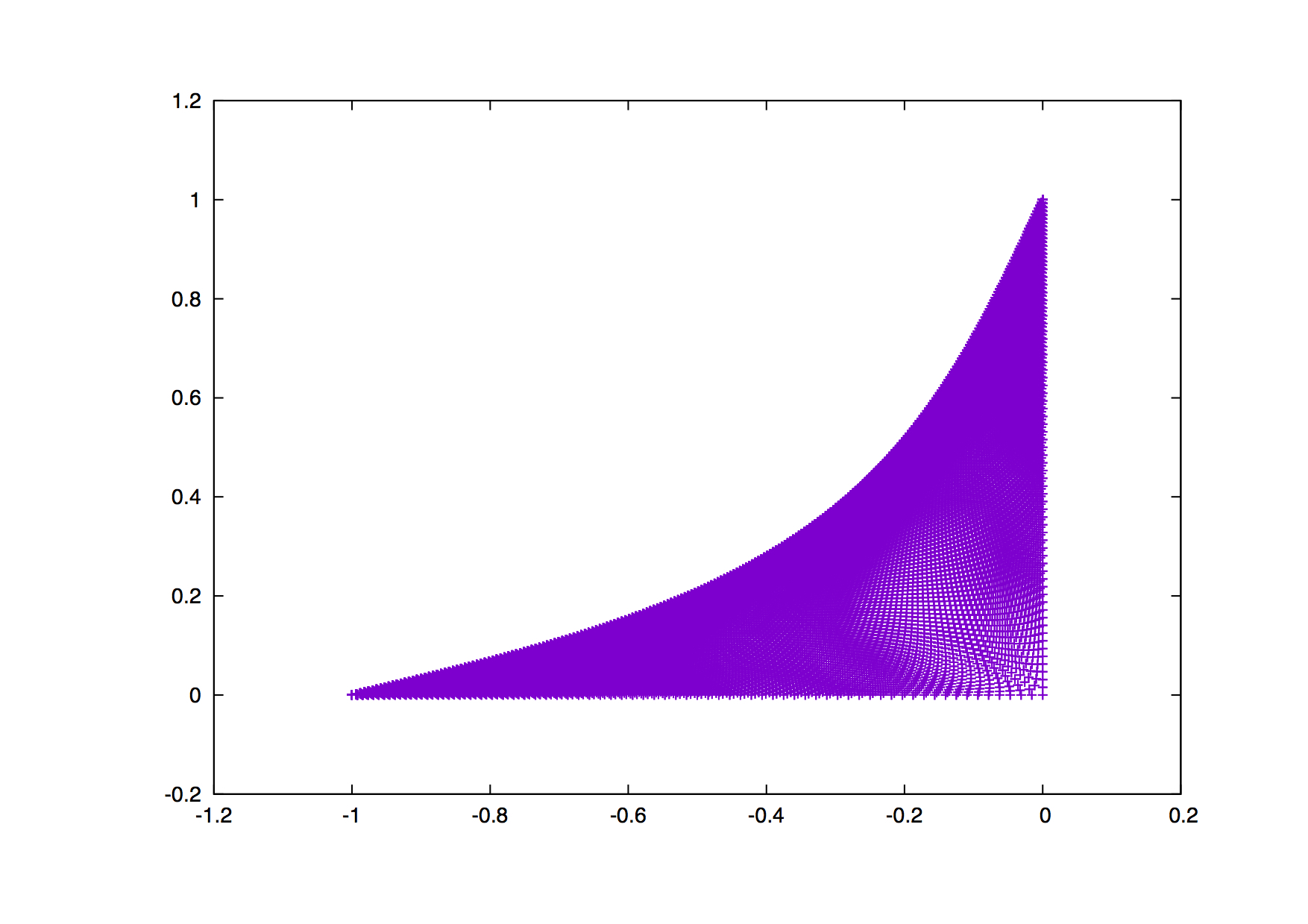}
    }
    \hfill
    \subfloat[{\label{fig:r1ia2}
        The image of the vectorial $8$-Laplacian.
    }]{
      \includegraphics[scale=\figscale,width=0.35\figwidth]{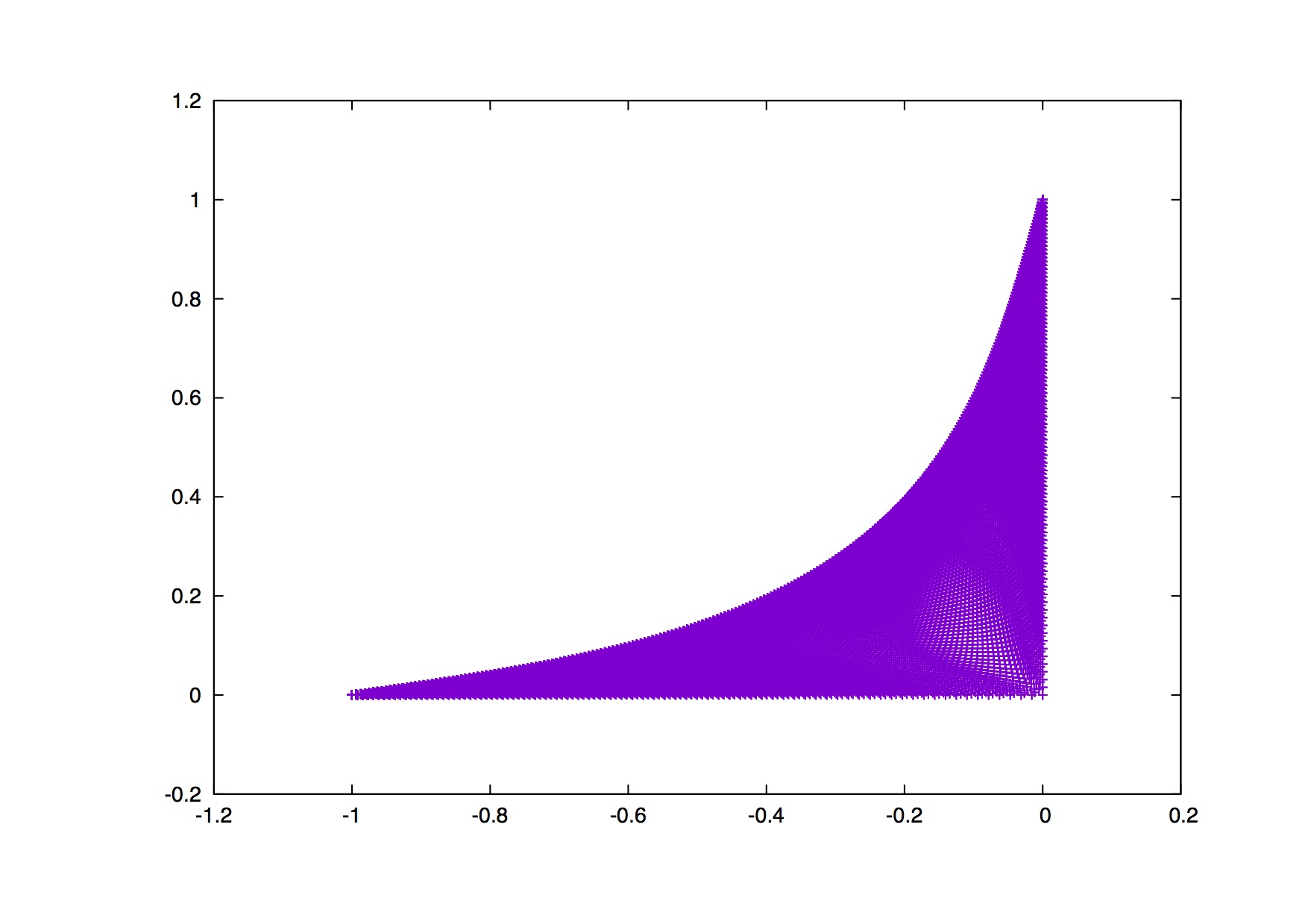}
    }
    \hfill
    \subfloat[{\label{fig:r1ia3}
        The image of the vectorial $64$-Laplacian.
    }]{
      \includegraphics[scale=\figscale,width=0.35\figwidth]{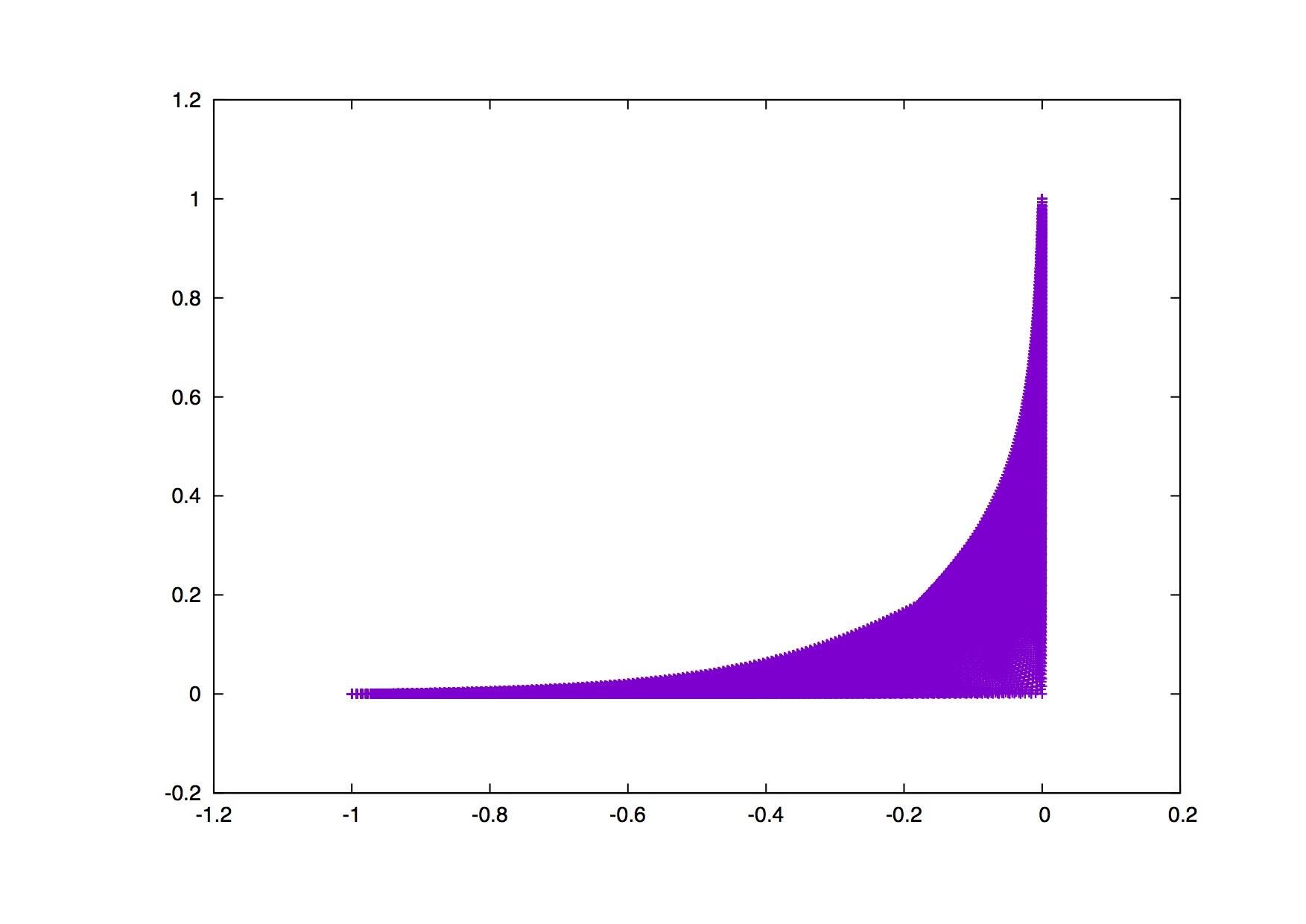}
    }
    \hfill
    \subfloat[{\label{fig:r1ia4}
        The image of the vectorial $256$-Laplacian.
    }]{
      \includegraphics[scale=\figscale,width=0.35\figwidth]{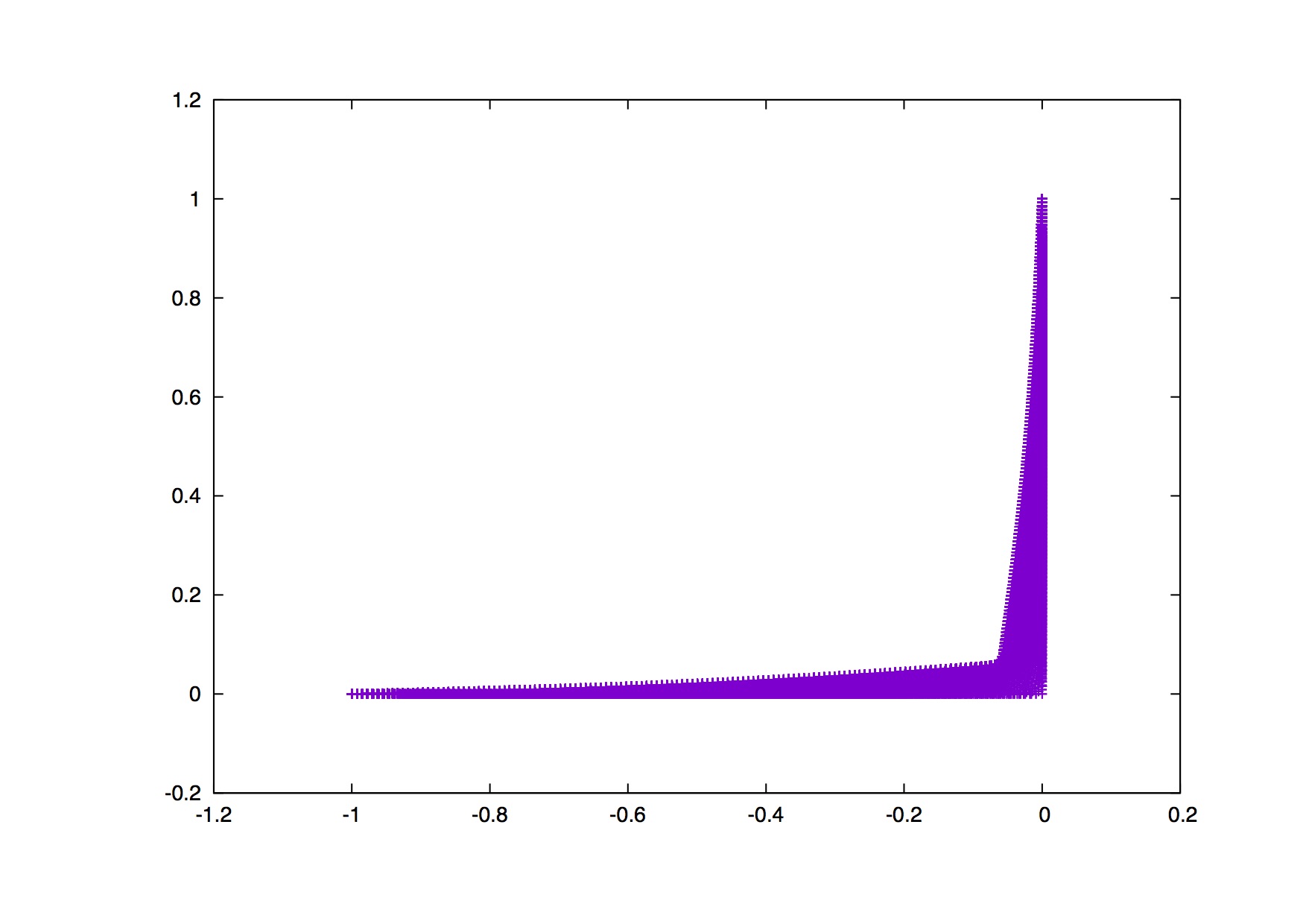}
    }
    \hfill
    \subfloat[{\label{fig:r1ia5}
        The image of the vectorial $512$-Laplacian.
    }]{
      \includegraphics[scale=\figscale,width=0.35\figwidth]{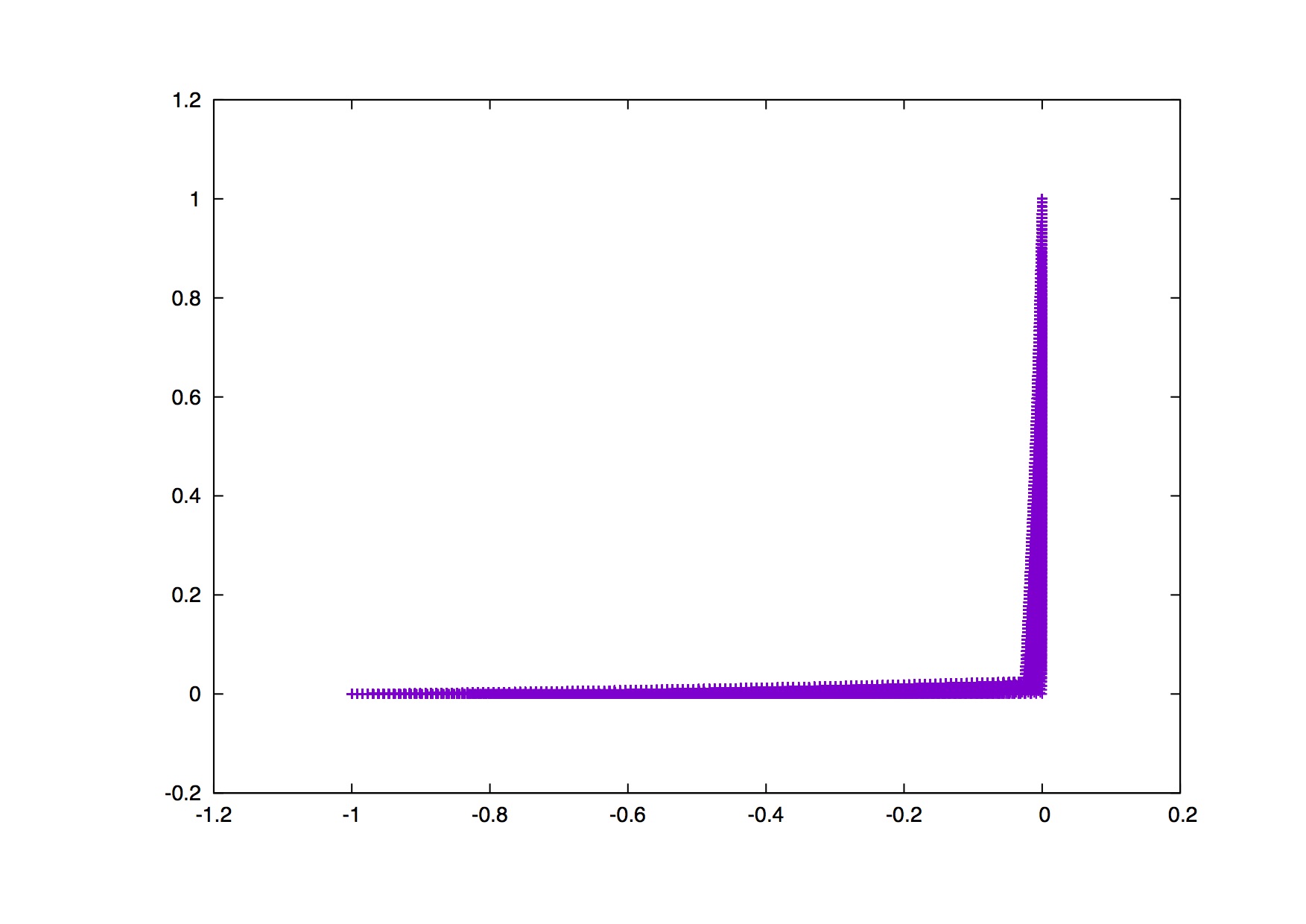}
    }
    \hfill
    \subfloat[{\label{fig:r1ia6}
        The image of the vectorial $1024$-Laplacian.
    }]{
      \includegraphics[scale=\figscale,width=0.35\figwidth]{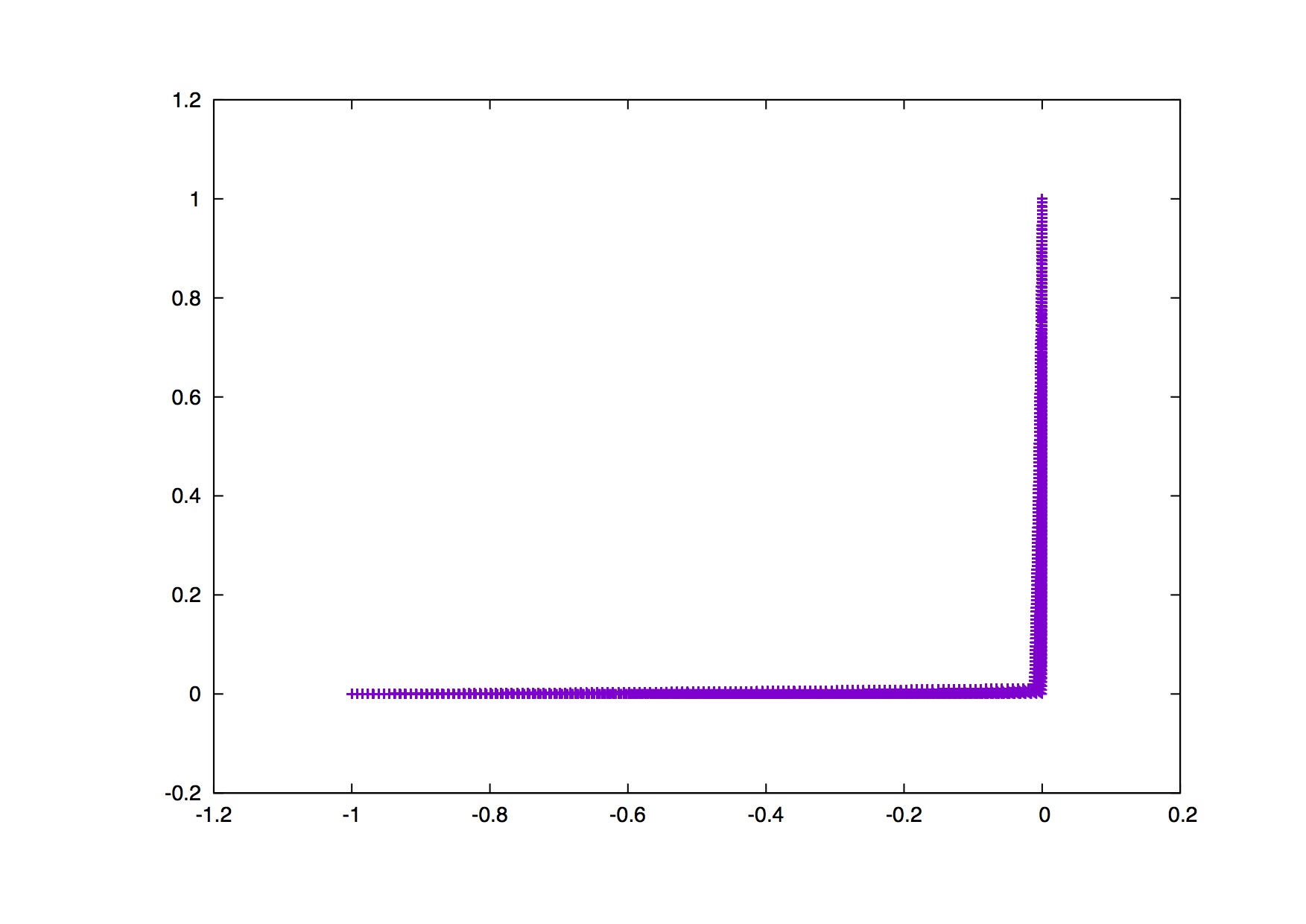}
    }
    \end{center}
  \end{figure}

\begin{figure}[!ht]
  \caption[Numerical Results for Problem \eqref{eqn:Problem:1} with
  $P^1$ elements]
  {\label{fig:rank1rk}
    An illustration of the rank of the solution to the vectorial $p$-Laplacian with the rank-one boundary conditions given in (\ref{eq:rank1-bcs}) for various increasing values of $p$. Here we plot $\det{\D U}$ and associated contour lines. These are plotted at increments of $0.05$. Notice as $p$ increases, the solution has lower rank over a larger portion of the domain, that is the size of $\W_1$ increases.
  }
  \begin{center}
    \subfloat[{\label{fig:rka1}
        $\det{\D U}$ for the vectorial $2$-Laplacian.
    }]{
      \includegraphics[scale=\figscale,width=0.35\figwidth]{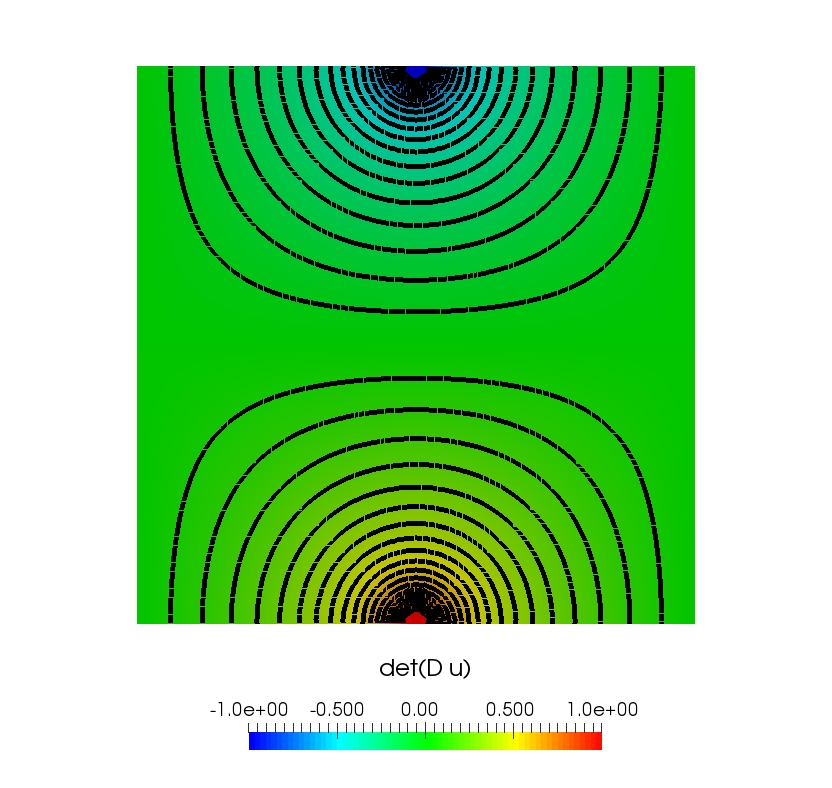}
    }
    \hfill
    \subfloat[{\label{fig:rka2}
        $\det{\D U}$ for the vectorial $8$-Laplacian.
    }]{
      \includegraphics[scale=\figscale,width=0.35\figwidth]{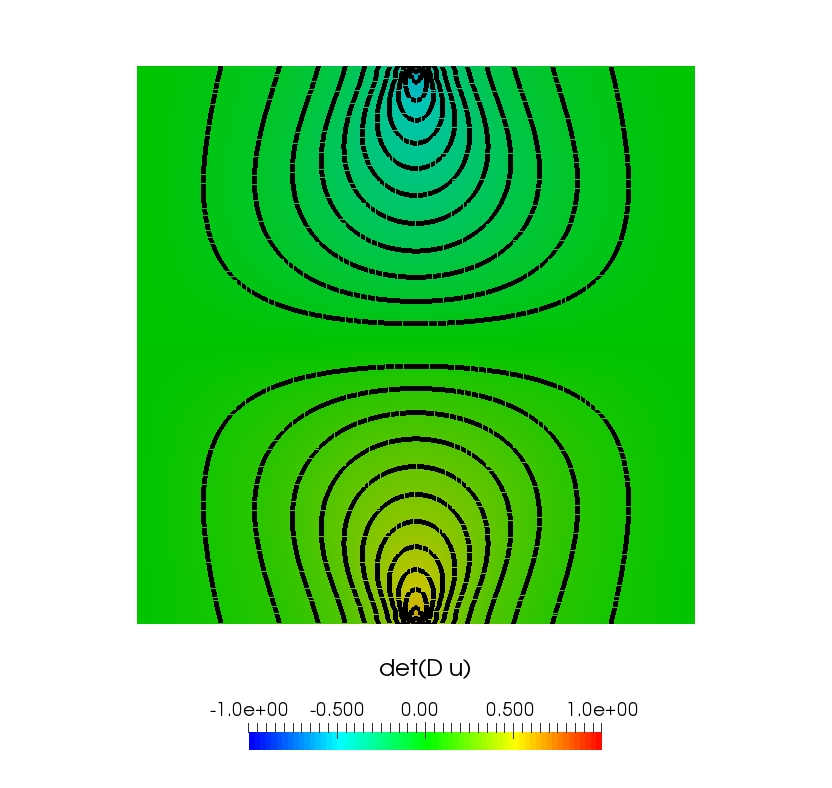}
    }
    \hfill
    \subfloat[{\label{fig:rka3}
        $\det{\D U}$ for the vectorial $64$-Laplacian.
    }]{
      \includegraphics[scale=\figscale,width=0.35\figwidth]{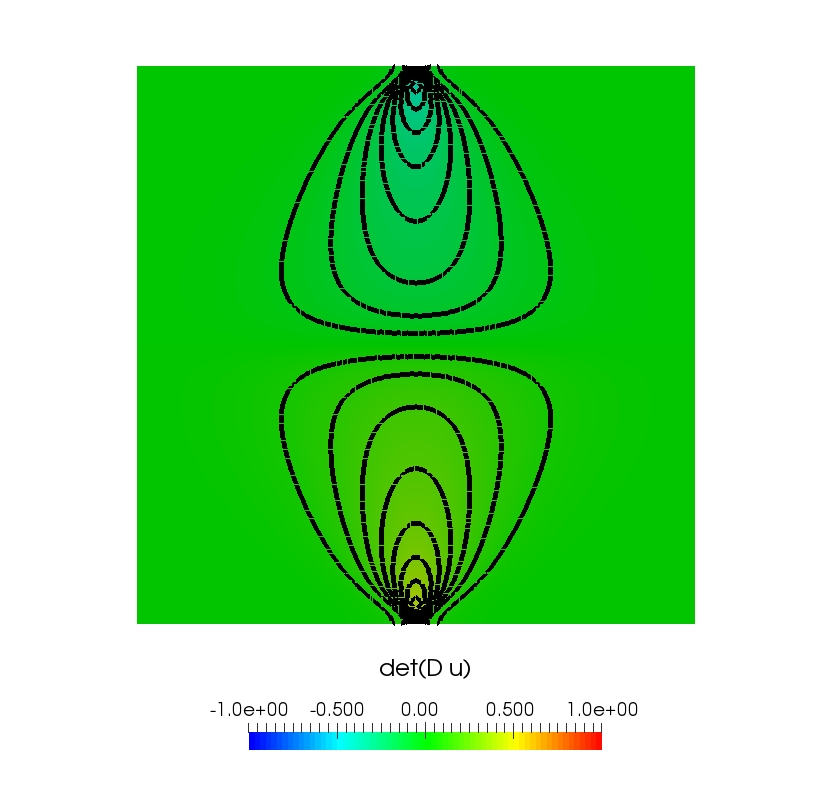}
    }
    \hfill
    \subfloat[{\label{fig:rka4}
        $\det{\D U}$ for the vectorial $256$-Laplacian.
    }]{
      \includegraphics[scale=\figscale,width=0.35\figwidth]{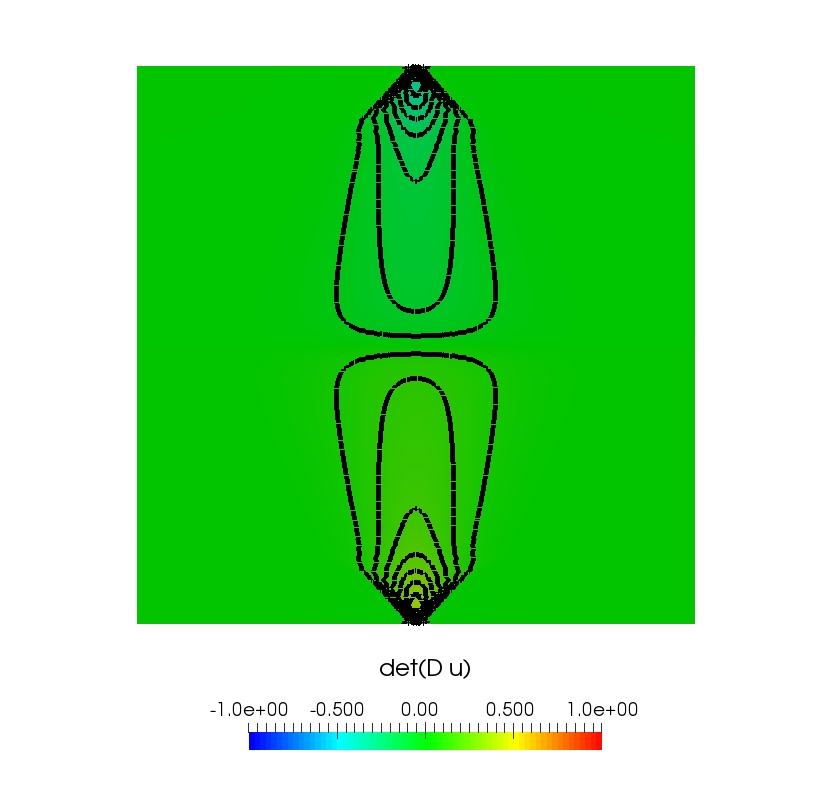}
    }
    \hfill
    \subfloat[{\label{fig:rka5}
        $\det{\D U}$ for the vectorial $512$-Laplacian.
    }]{
      \includegraphics[scale=\figscale,width=0.35\figwidth]{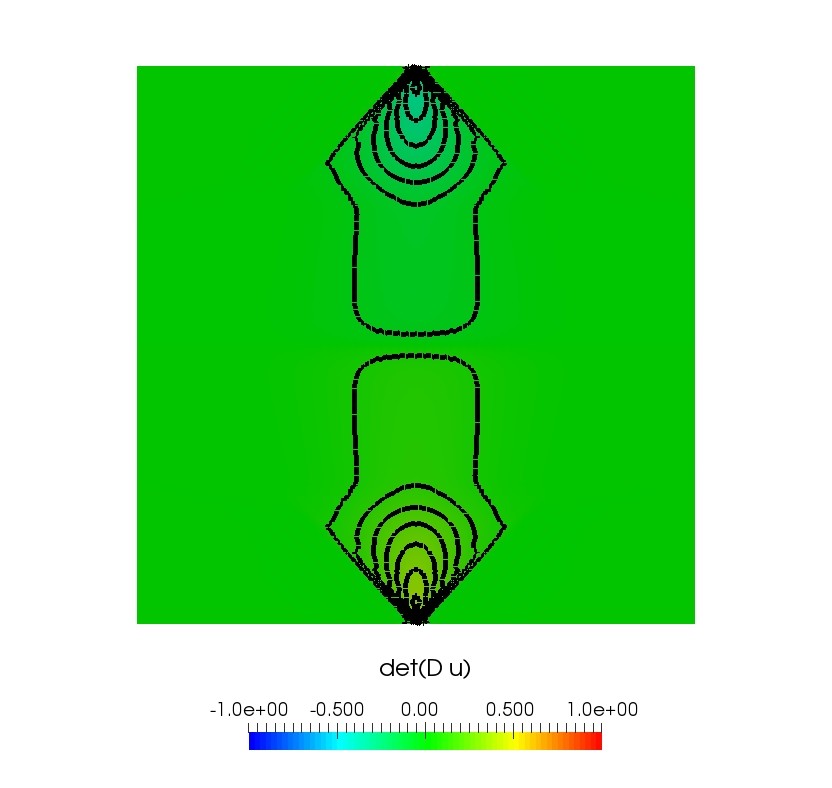}
    }
    \hfill
    \subfloat[{\label{fig:rka6}
        $\det{\D U}$ for the vectorial $1024$-Laplacian.
    }]{
      \includegraphics[scale=\figscale,width=0.35\figwidth]{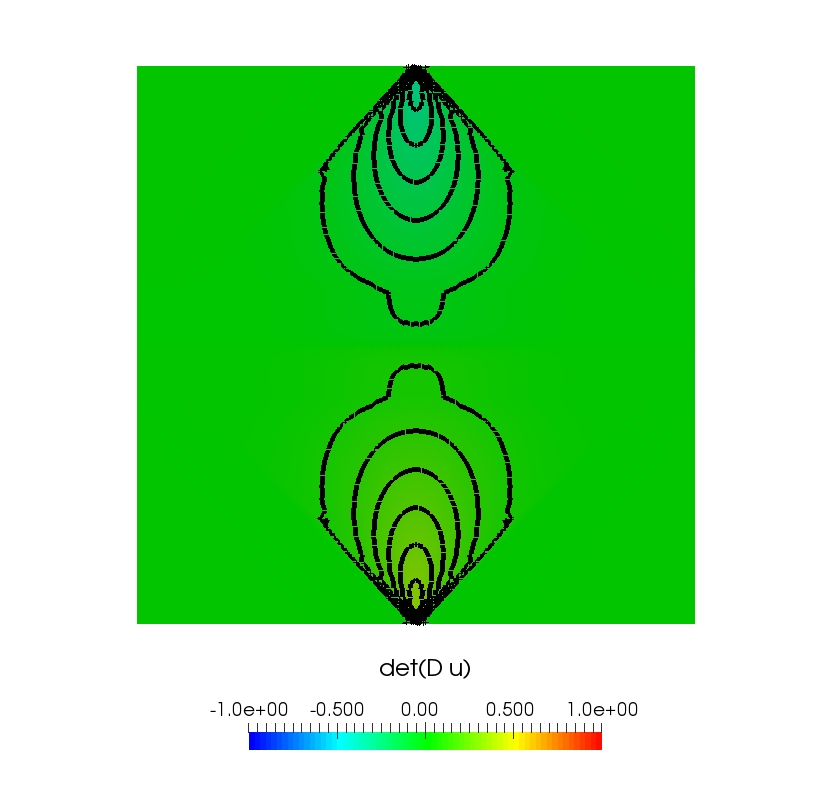}
    }
    \end{center}
  \end{figure}

In this test we construct mappings with rank-one boundary data. We take
\begin{equation}
  \label{eq:rank1-bcs}
  g(x,y) \, :=\, 
  \begin{cases}
    \vec e_x x,\ \ \text{ if } x < 0, \ms
    \\
    \vec e_y x, \ \ \text{ if } x \geq 0,
  \end{cases}
\end{equation}
  where $\vec e_x = (1,0)^\top$ and $\vec e_y = (0,1)^\top$ denote the unit vectors in the $x$ and $y$ directions. Notice that $g \in \cont{0}(\partial\W,\R^2)$. We examine the image and rank of the numerical approximation for various values of $p$ in Figures \ref{fig:rank1im} and \ref{fig:rank1rk}. 

\subsection{Solutions $(-1,+1)^2\sub \R^2\larrow  \R^2$ with boundary data an explicit $\infty$-Harmonic map with triple junction interface.} \label{subsection3.5}

\begin{figure}[!h]
  \caption[]
  {\label{fig:trip}
    An illustration of the rank of the solution to the vectorial $p$-Laplacian with the rank-one boundary conditions given in (\ref{eq:rank1-bcs}) for various increasing values of $p$. Here we plot $\det{\D U}$ and associated contour lines. These are plotted at increments of $0.05$. Notice as $p$ increases, the structure illustrated in Figure \ref{fig:1a} becomes more pronounced. 
  }
  \begin{center}
    \subfloat[{\label{fig:trip1}
        $\det{\D U}$ for the vectorial $2$-Laplacian.
    }]{
      \includegraphics[scale=\figscale,width=0.35\figwidth]{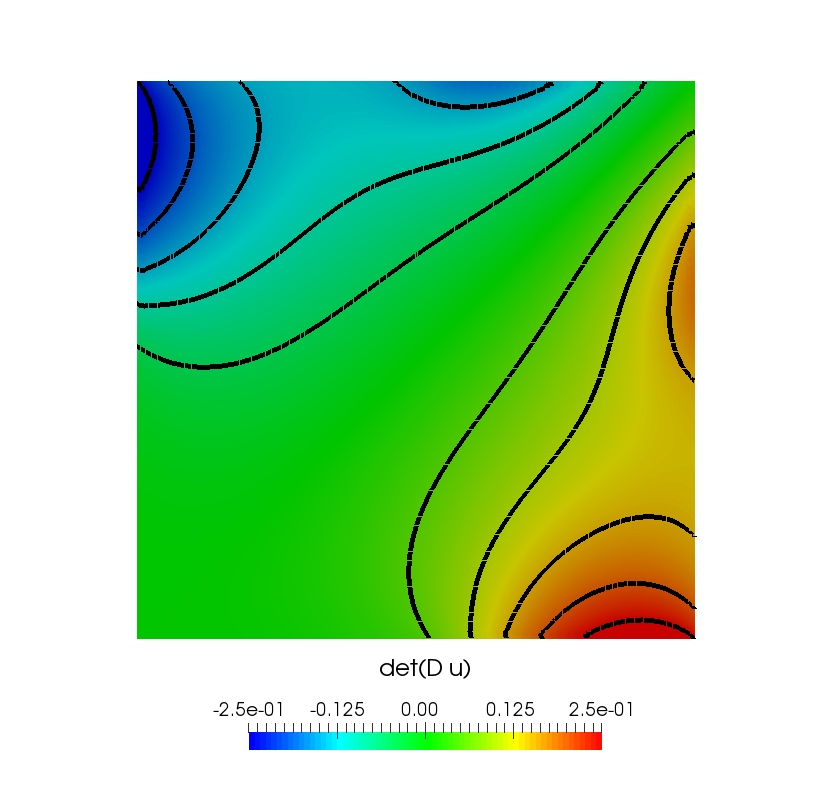}
    }
    \hfill
    \subfloat[{\label{fig:trip2}
        $\det{\D U}$ for the vectorial $8$-Laplacian.
    }]{
      \includegraphics[scale=\figscale,width=0.35\figwidth]{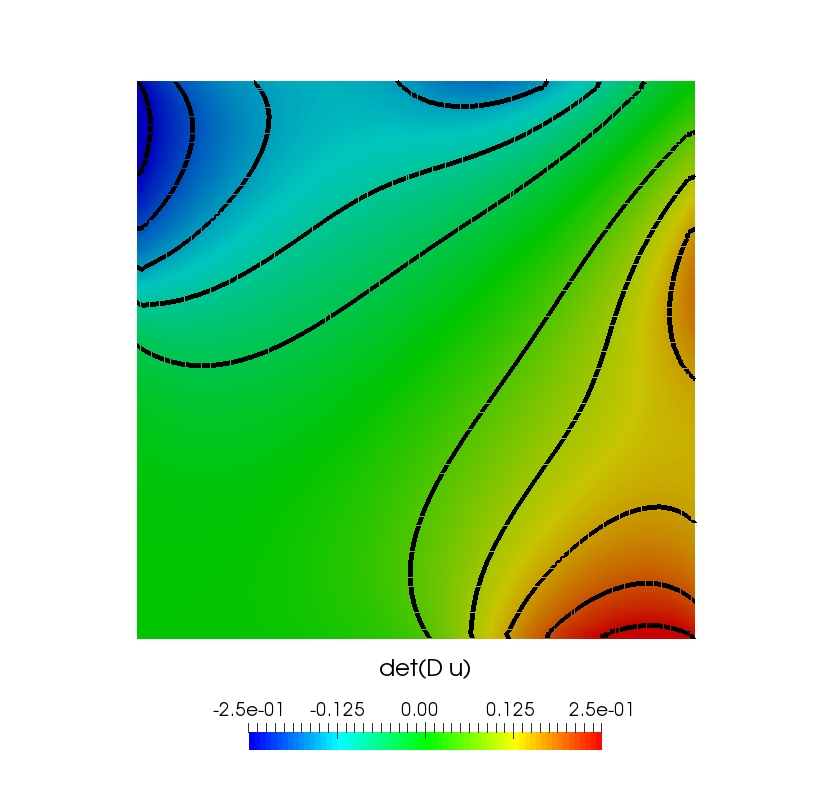}
    }
    \hfill
    \subfloat[{\label{fig:trip3}
        $\det{\D U}$ for the vectorial $64$-Laplacian.
    }]{
      \includegraphics[scale=\figscale,width=0.35\figwidth]{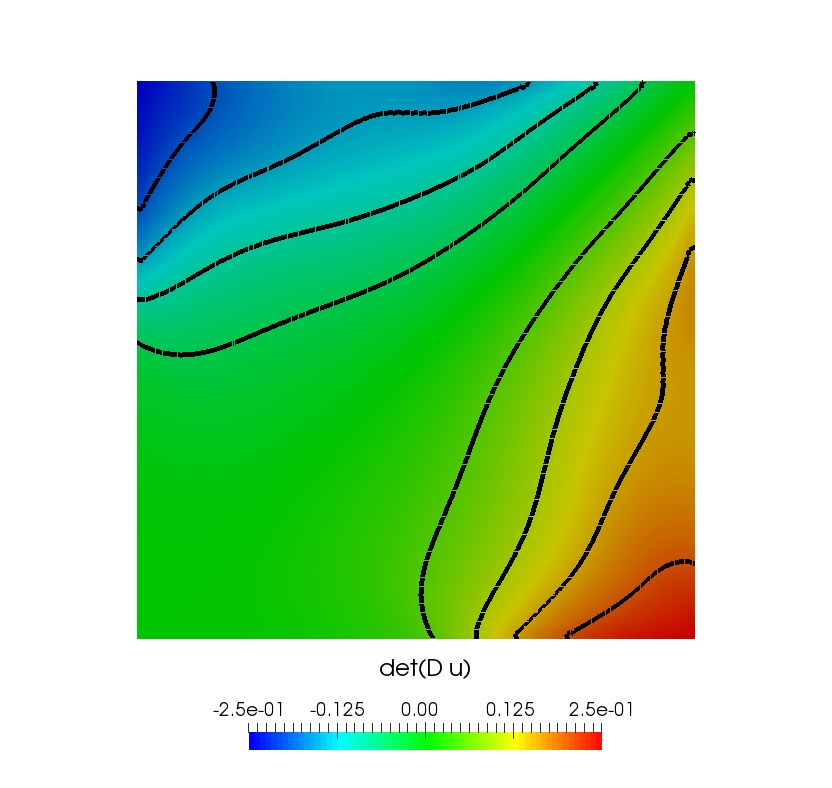}
    }
    \hfill
    \subfloat[{\label{fig:trip4}
        $\det{\D U}$ for the vectorial $256$-Laplacian.
    }]{
      \includegraphics[scale=\figscale,width=0.35\figwidth]{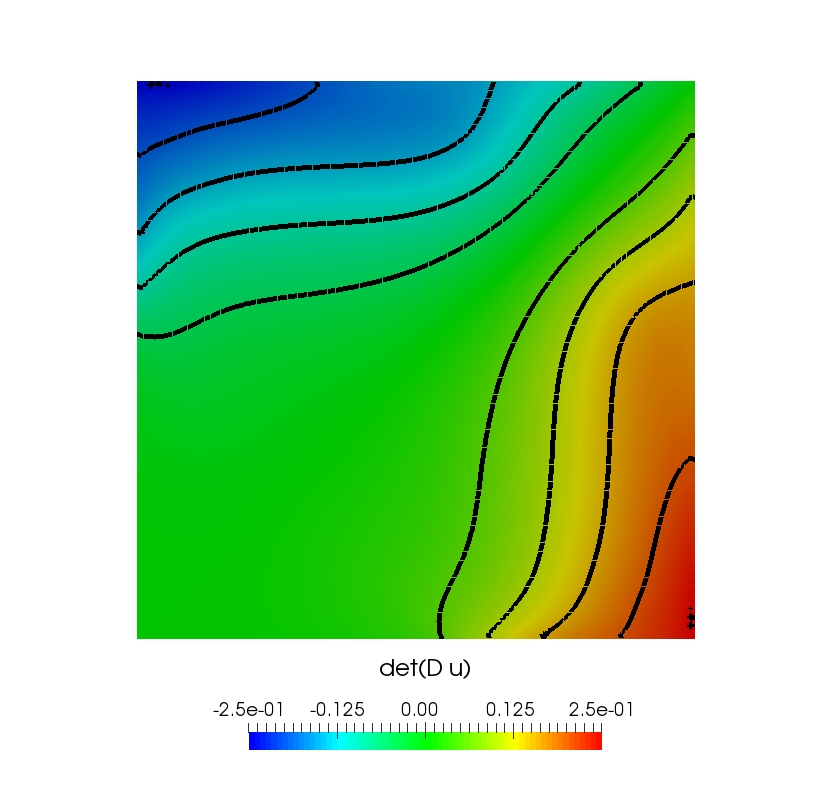}
    }
    \hfill
    \subfloat[{\label{fig:trip5}
        $\det{\D U}$ for the vectorial $512$-Laplacian.
    }]{
      \includegraphics[scale=\figscale,width=0.35\figwidth]{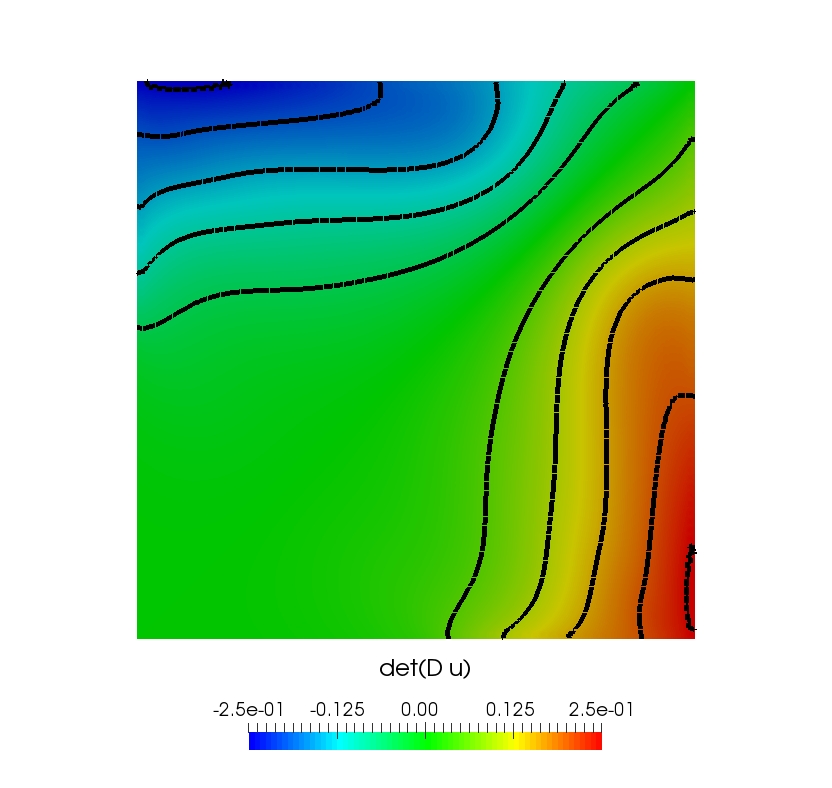}
    }
    \hfill
    \subfloat[{\label{fig:trip6}
        $\det{\D U}$ for the vectorial $1024$-Laplacian.
    }]{
      \includegraphics[scale=\figscale,width=0.35\figwidth]{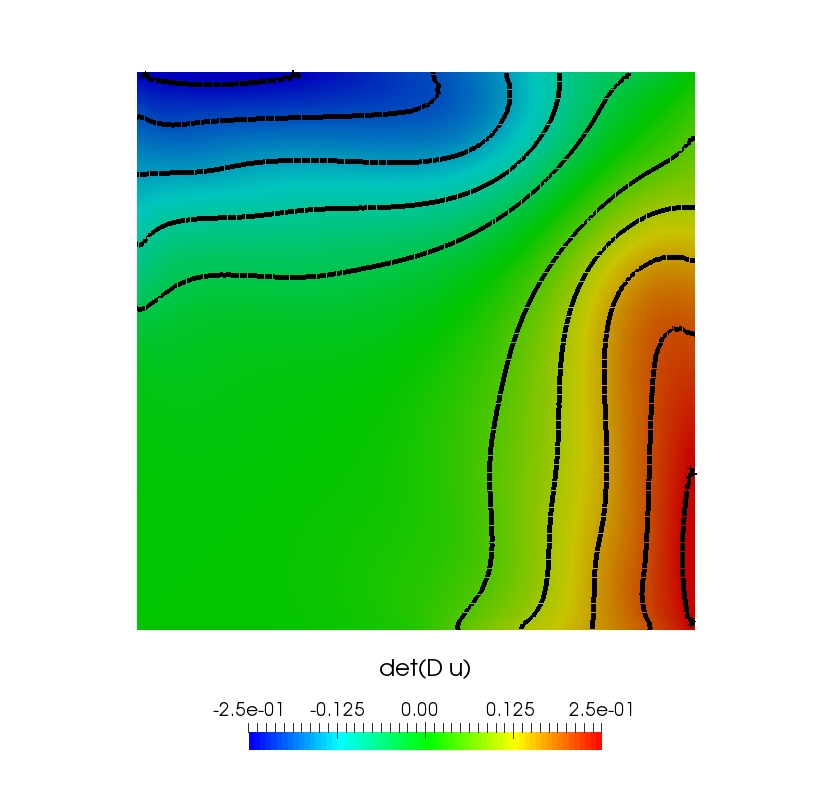}
    }
    \end{center}
  \end{figure}

In this test we construct boundary data which give rise to the example of an explicit smooth $\infty$-Harmonic mapping with a triple junction interface as illustrated in Figure \ref{fig:1a}. We take
\begin{equation}
  K(t) \,:=\,
  \begin{cases}
    1-\dfrac{1}{1+t^2}, \text{ if } t > 0,
    \\
    0, \ \ \ \qquad \ \ \ \text{ otherwise.}
  \end{cases}
\end{equation}
We also take
\begin{equation}
  g(x,y) \,:=\,
  \frac{3}{4}\qp{\int_y^x \cos{K(t)}\text{d}t ,\int_y^x \sin{K(t)}\text{d}t }^{\!\!\top}.
\end{equation}
The numerical experiment is given in Figure \ref{fig:trip}.

\subsection{Solutions $(-1,+1)^2\sub \R^2\larrow  \R^2$ with boundary data an explicit $\infty$-Harmonic map with rectangular interface.} \label{subsection3.4}

\begin{figure}[!h]
  \caption[Numerical Results for Problem \eqref{eqn:Problem:1} with
  $P^1$ elements]
  {\label{fig:box}
    An illustration of the rank of the solution to the vectorial $p$-Laplacian with the boundary conditions given in \ref{2.8} for the rectangular interface illustrated in Figure \ref{fig:1a} for various increasing values of $p$. Here we plot $\det{\D U}$ and associated contour lines. These are plotted at increments of $0.05$. Notice as $p$ increases, the structure illustrated in Figure \ref{fig:1b} becomes more pronounced. 
  }
  \begin{center}
    \subfloat[{\label{fig:a1}
        $\det{\D U}$ for the vectorial $2$-Laplacian.
    }]{
      \includegraphics[scale=\figscale,width=0.35\figwidth]{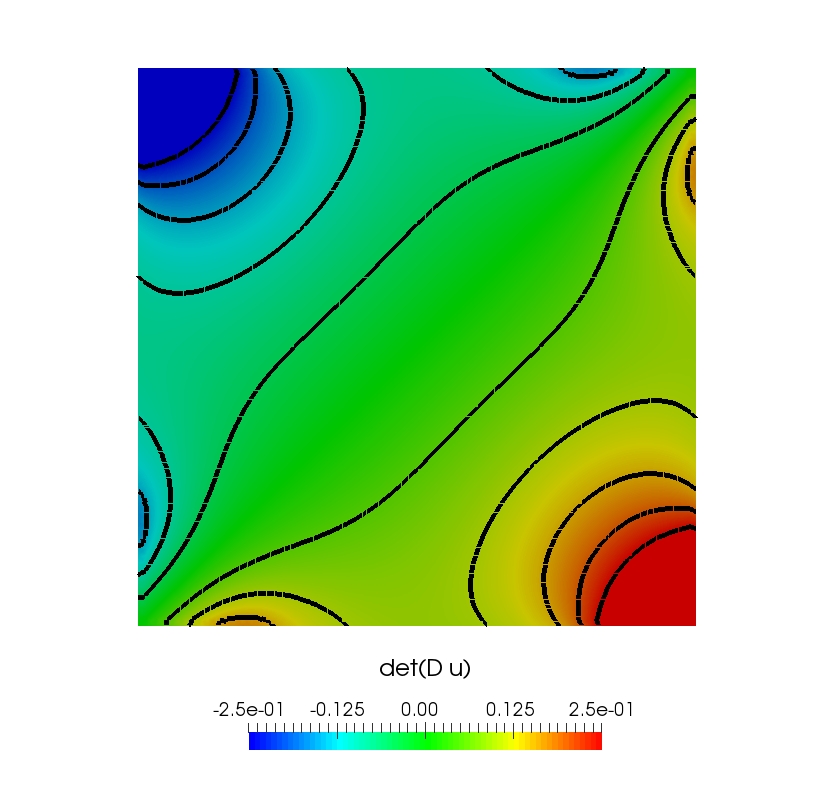}
    }
    \hfill
    \subfloat[{\label{fig:a2}
        $\det{\D U}$ for the vectorial $8$-Laplacian.
    }]{
      \includegraphics[scale=\figscale,width=0.35\figwidth]{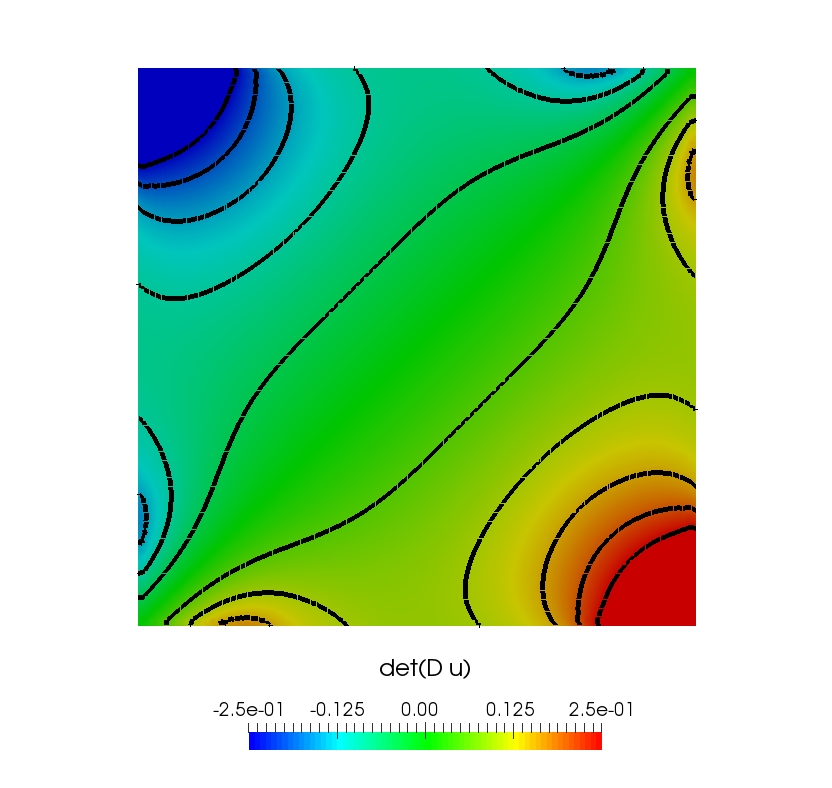}
    }
    \hfill
    \subfloat[{\label{fig:a2}
        $\det{\D U}$ for the vectorial $64$-Laplacian.
    }]{
      \includegraphics[scale=\figscale,width=0.35\figwidth]{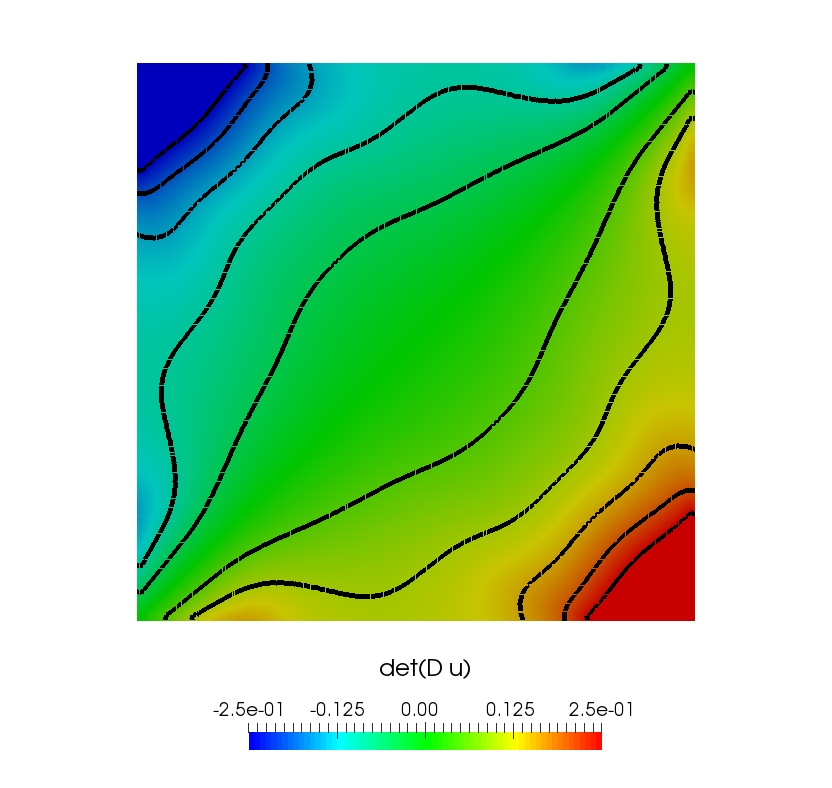}
    }
    \hfill
    \subfloat[{\label{fig:a2}
        $\det{\D U}$ for the vectorial $256$-Laplacian.
    }]{
      \includegraphics[scale=\figscale,width=0.35\figwidth]{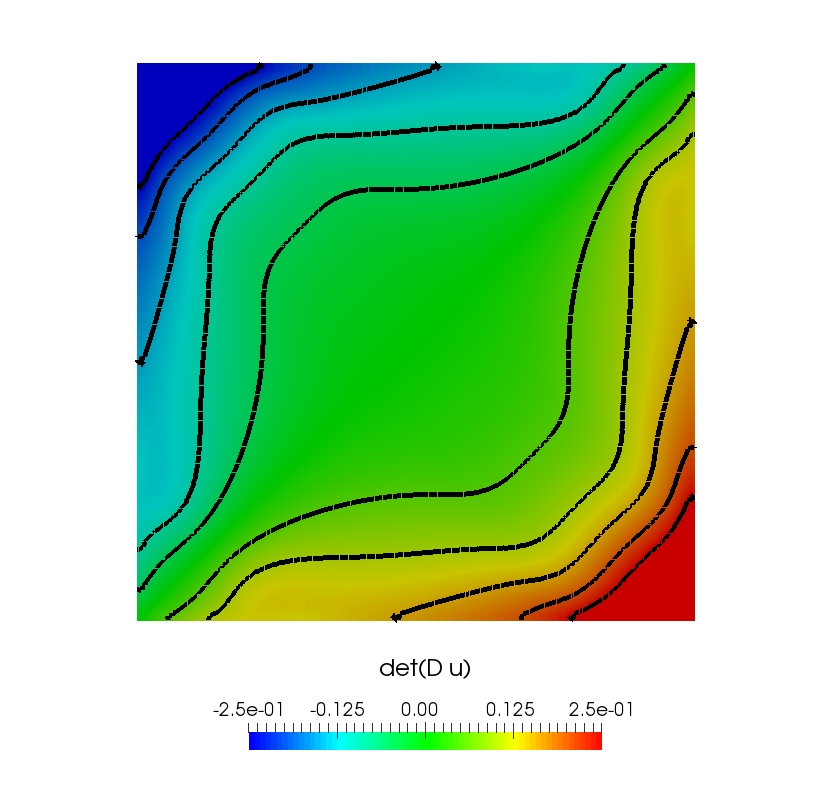}
    }
    \hfill
    \subfloat[{\label{fig:a2}
        $\det{\D U}$ for the vectorial $512$-Laplacian.
    }]{
      \includegraphics[scale=\figscale,width=0.35\figwidth]{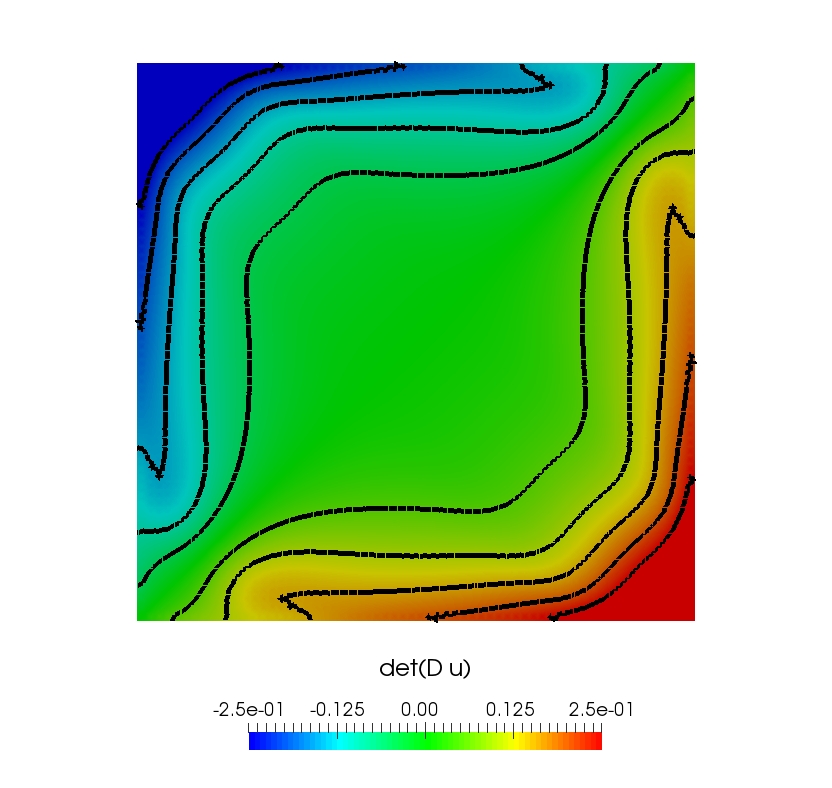}
    }
    \hfill
    \subfloat[{\label{fig:a2}
        $\det{\D U}$ for the vectorial $1024$-Laplacian.
    }]{
      \includegraphics[scale=\figscale,width=0.35\figwidth]{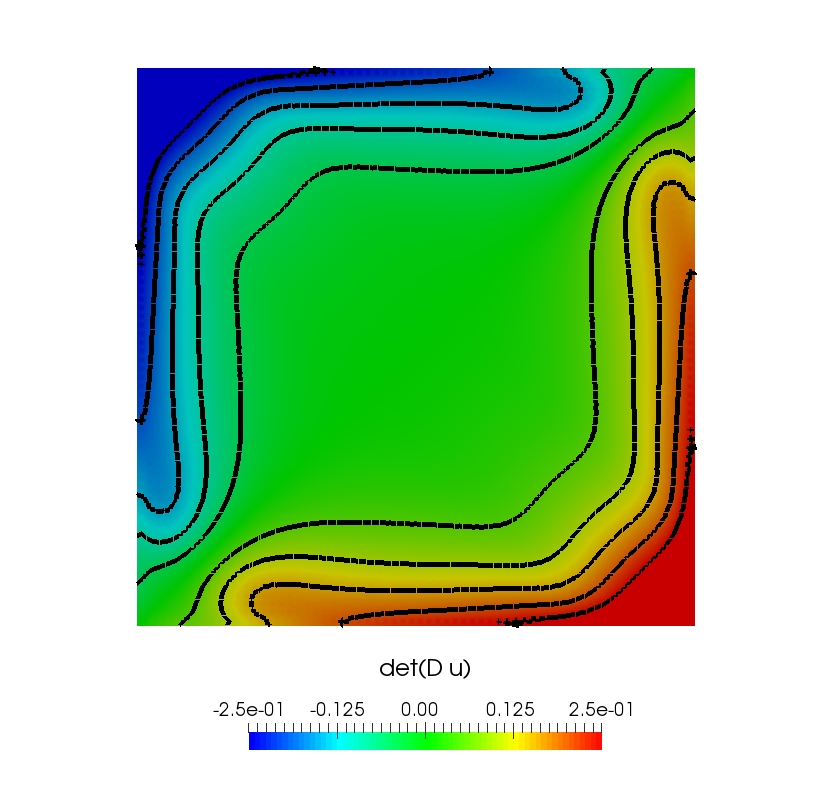}
    }
    \end{center}
  \end{figure}

In this test we constuct boundary data which give rise to the example of an explicit smooth $\infty$-Harmonic mapping with a box interface as illustrated in Figure \ref{fig:1b}.
We take
\begin{equation}
  K(t) \,:=\,
\left\{
  \begin{array}{ll}
    1-\dfrac{1}{1+\qp{t-1}^2+1}, & \text{ if } t > 1, \smallskip
    \\
    \dfrac{1}{1+\qp{t+1}^2+1} -1, & \text{ if } t < -1, \smallskip 
    \\
    0, & \text{ otherwise.}
  \end{array}
  \right.
\end{equation}
and
\begin{equation}
  g(x,y) \, :=\,
  \frac{3}{4} \qp{\int_y^x \cos{K(t)}\text{d}t,\int_y^x \sin{K(t)} \text{d}t}^{\!\!\top}.
\end{equation}
The numerical experiment is given in Figure \ref{fig:box}.

\clearpage


\begin{thebibliography}{10}

\bibitem[A1]{A1} G. Aronsson, \emph{Minimization problems for the functional $sup_x F(x,
f(x), f'(x))$}, Arkiv f\"ur Mat. 6 (1965), 33 - 53.

\bibitem[A2]{A2} G. Aronsson, \emph{Minimization problems for the functional $sup_x F(x,
f(x), f'(x))$ II}, Arkiv f\"ur Mat. 6 (1966), 409 - 431.

\bibitem[A3]{A3} G. Aronsson, \emph{Extension of functions satisfying Lipschitz conditions}, Arkiv f\"ur Mat. 6 (1967), 551 - 561.

\bibitem[A4]{A4} G. Aronsson, \emph{On the partial differential equation $u_x^2 u_{xx} + 2u_x u_y u_{xy} + u_y^2 u_{yy} = 0$}, Arkiv f\"ur Mat. 7
(1968), 395 - 425.

\bibitem[A5]{A5} G. Aronsson, \emph{Minimization problems for the functional $sup_x F(x,
f(x), f'(x))$ III}, Arkiv f\"ur Mat. (1969), 509 - 512.

\bibitem[A6]{A6} G. Aronsson, \emph{On Certain Singular Solutions of the Partial
Differential Equation $u_x^2 u_{xx} + 2u_x u_y u_{xy} + u_y^2 u_{yy} =
0$}, Manuscripta Math. 47 (1984), no 1-3, 133 - 151.

\bibitem[A7]{A7} G. Aronsson, \emph{Construction of Singular Solutions to the $p$-Harmonic
Equation and its Limit Equation for $p=\infty$}, Manuscripta Math.
56 (1986), 135 - 158.

\bibitem[AS]{ArmstrongSmart:2010}
S.~Armstrong and C.~Smart.
\newblock An easy proof of jensen's theorem on the uniqueness of infinity
  harmonic functions.
\newblock {\em Calculus of Variations and Partial Differential Equations},
  37(3-4):381--384, 2010.

\bibitem[B]{B} N. Barron, \emph{Viscosity Solutions and Analysis in $L^\infty$},  Nonlinear analysis, differential equations and control (Montreal QC, 1998), Kluwer Acad. Publ. Dordrecht, 1999, 1 - 60.

\bibitem[BL]{BL} J. W. Barrett and W. B. Liu, \emph{Finite element approximation of some degenerate monotone quasilinear elliptic systems}, SIAM journal of numerical analysis 33 (1996), no 1, 88 - 106.

\bibitem[BS]{BarlesSouganidis:1991}
G.~Barles and P.~E. Souganidis.
\newblock Convergence of approximation schemes for fully nonlinear second order
  equations.
\newblock {\em Asymptotic Anal.}, 4(3):271--283, 1991.
  
\bibitem[C]{C} M. G. Crandall, \emph{A visit with the $\infty$-Laplacian}, in \emph{Calculus of Variations and Non-Linear PDE}, Springer Lecture notes in Mathematics 1927, CIME, Cetraro Italy 2005.

\bibitem[CIL]{CIL} M. G. Crandall, H. Ishii, P.-L. Lions, \emph{User's Guide to
Viscosity Solutions of 2nd Order Partial Differential Equations},
Bulletin of the AMS, Vol. 27, Nr 1, Pages 1 - 67, 1992.

\bibitem[CT]{CrouzeixThomee:1987}
M.~Crouzeix and V.~Thom{\'e}e, \emph{The stability in {$L\sb p$} and {$W\sp
  1\sb p$} of the {$L\sb 2$}-projection onto finite element function spaces},
  Math. Comp. {48} (1987), 521--532. 

\bibitem[D]{Dacorogna:2008}
B. Dacorogna.
\newblock {\em Direct methods in the calculus of variations}, volume~78 of {\em
  Applied Mathematical Sciences}.
\newblock Springer, New York, second edition, 2008.

\bibitem[EO]{EsedogluOberman:2011}
S. Esedoglu and A. Oberman.
\newblock Fast semi-implicit solvers for the infinity laplace and p-laplace
  equations.
\newblock {\em Arxiv}, 2011.

\bibitem[HLL]{HuangLiLiu:2007}
Y.~Q. Huang, R. Li, and W. Liu, \emph{Preconditioned descent algorithms
  for {$p$}-{L}aplacian}, J. Sci. Comput. {32} (2007), 343--371.

\bibitem[J]{Jensen:1993}
R. Jensen.
\newblock Uniqueness of {L}ipschitz extensions: minimizing the sup norm of the
  gradient.
\newblock {\em Arch. Rational Mech. Anal.}, 123(1):51--74, 1993.

\bibitem[K]{K} N. Katzourakis, \emph{An Introduction to viscosity Solutions for Fully Nonlinear PDE with Applications to Calculus of Variations in $L^\infty$}, Springer Briefs in Mathematics, 2015, DOI 10.1007/978-3-319-12829-0.

\bibitem[K1]{K1} \_\_\_\_\_\_\_\_\_\_\_\_ ,  \emph{$L^\infty$-Variational Problems for Maps and the Aronsson PDE system}, J.\ Differential Equations, Volume 253, Issue 7 (2012), 2123 - 2139.

\bibitem[K2]{K2} \_\_\_\_\_\_\_\_\_\_\_\_ ,  \emph{Explicit $2D$ $\infty$-Harmonic Maps whose Interfaces have Junctions and Corners}, Comptes Rendus Acad. Sci. Paris, Ser.I, 351 (2013) 677 - 680.

\bibitem[K3]{K3} \_\_\_\_\_\_\_\_\_\_\_\_ ,   \emph{On the Structure of $\infty$-Harmonic Maps}, Communications in PDE, Volume 39, Issue 11 (2014), 2091 - 2124.

\bibitem[K4]{K4} \_\_\_\_\_\_\_\_\_\_\_\_ , \emph{$\infty$-Minimal Submanifolds}, Proceedings of the Amer. Math. Soc., 142 (2014) 2797-2811.

\bibitem[K5]{K5} \_\_\_\_\_\_\_\_\_\_\_\_ ,   \emph{Nonuniqueness in Vector-valued Calculus of Variations in $L^\infty$ and some Linear Elliptic Systems}, Communications on Pure and Applied Analysis,  Vol. 14, 1, 313 - 327 (2015). 

\bibitem[K6]{K6} \_\_\_\_\_\_\_\_\_\_\_\_ ,    \emph{Optimal $\infty$-Quasiconformal Immersions}, ESAIM Control, Opt. and Calc. Var., to appear (2015) DOI: http://dx.doi.org/10.1051/cocv/2014038. 

\bibitem[K7]{K7} \_\_\_\_\_\_\_\_\_\_\_\_ ,  \emph{Generalised solutions for fully nonlinear PDE systems and existence-uniqueness theorems}, ArXiv preprint, \url{http://arxiv.org/pdf/1501.06164.pdf}. 

\bibitem[K8]{K8} \_\_\_\_\_\_\_\_\_\_\_\_ ,   \emph{Existence of generalised solutions to the equations of vectorial Calculus of Variations in $L^\infty$}, ArXiv preprint, \url{http://arxiv.org/pdf/1502.01179.pdf}. 

\bibitem[K9]{K9} \_\_\_\_\_\_\_\_\_\_\_\_ ,   \emph{A New Characterisation of $\infty$-Harmonic and $p$-Harmonic Mappings via Affine Variations in $L^\infty$}, ArXiv preprint, \url{http://arxiv.org/pdf/1509.01811.pdf}.

\bibitem[K10]{K10} \_\_\_\_\_\_\_\_\_\_\_\_ ,  \emph{Equivalence between weak and $\mD$-solutions for symmetric hyperbolic PDE systems}, ArXiv preprint, \url{http://arxiv.org/pdf/1507.03042.pdf}. 

\bibitem[K11]{K11} \_\_\_\_\_\_\_\_\_\_\_\_ ,  \emph{Mollification of $\mD$-solutions to fully nonlinear PDE systems}, ArXiv preprint, \url{http://arxiv.org/pdf/1508.05519.pdf}. 

\bibitem[LP]{LakkisPryer:2013}
O.~Lakkis and T.~Pryer.
\newblock A finite element method for nonlinear elliptic problems.
\newblock {\em SIAM Journal on Scientific Computing}, 35(4):A2025--A2045, 2013.

\bibitem[LP1]{Pryer:2013}
O. Lakkis and T. Pryer, \emph{An adaptive finite element method for the infinity laplacian}, Numerical Mathematics and Advanced Applications, 283 - 291 (2013).


\bibitem[O]{Oberman:2005}
A. Oberman.
\newblock A convergent difference scheme for the infinity {L}aplacian:
  construction of absolutely minimizing {L}ipschitz extensions.
\newblock {\em Math. Comp.}, 74(251):1217--1230 (electronic), 2005.

\bibitem[O1]{Oberman:2013}
A. Oberman.
\newblock Finite difference methods for the infinity {L}aplace and
  {$p$}-{L}aplace equations.
\newblock {\em J. Comput. Appl. Math.}, 254:65--80, 2013.

\bibitem[P]{P} T. Pryer, \emph{On the finite element approximation of Infinity-Harmonic functions}, ArXiv preprint, \url{http://arxiv.org/pdf/1511.00471}.
  
\bibitem[SS]{SS} S. Sheffield, C.K. Smart, \emph{Vector Valued Optimal Lipschitz Extensions}, Comm. Pure Appl. Math., Vol. 65, Issue 1, January 2012, 128 - 154.

\end{thebibliography}
\end{document}